\documentclass[11pt,epsfig]{article}
\usepackage{amsmath}
\usepackage{amsthm}
\usepackage{epsfig}
\usepackage{graphicx}
\usepackage{hyperref}
\usepackage{amsfonts}
\usepackage{graphicx, color, colortbl, epstopdf}

\newtheoremstyle{thmm}{1.5ex plus 1ex minus .2ex}{1.5ex plus 1ex minus .2ex}{\rmfamily}{}{\bfseries}{}{1em}{}
\theoremstyle{thmm}
\newtheorem{theorem}{Theorem}[section]
\newtheorem{lemma}{Lemma}[section]
\newtheorem{corollary}{Corollary}[section]

\renewenvironment{proof}[1][Proof]{\noindent\textit{#1. } }{\hfill$\square$}
\newcommand{\nn}{\nonumber}
\def \endproof{\vrule height8pt width 5pt depth 0pt}
\def\refe#1{(\ref{#1})}

\def\v{\varepsilon}
\def\e{\epsilon}

\def\u{{\bf u}} 
 
\def\wt{\widetilde}  
  
\def\u{{\bf u}}

\def\H{{\bf H}}
\def\v{{\bf v}} 
\def\C{{\cal C}}

\begin{document}

\date

\title{\bf Analysis of lowest-order characteristics-mixed FEMs for incompressible miscible flow in porous media}
\author{Weiwei Sun
\footnote{Advanced Institute of Natural Science, Beijing Normal University at Zhuhai, P.R. 
China  and Division of Science and Technology, United International 
College (BNU-HKBU), Zhuhai, 519087, P.R. China ({\it maweiw@uic.edu.cn}).
The work of this author was partially supported by a grant from National 
Natural Science Foundation of China under grant number 12071040, 
start-up funds (R5202009, R72021111) from United International College (BNU-HKBU) 
and Zhujiang Scholar program. 
} }
%

\maketitle

\begin{abstract}
The time discrete scheme of characteristics type is especially effective for convection-dominated diffusion problems. The scheme has been used in various engineering areas with different approximations  in spatial direction. The lowest-order mixed method is the most popular one for miscible flow in porous media. The method is based on a linear Lagrange approximation to the concentration and the zero-order Raviart-Thomas approximation to the pressure/velocity. However, the optimal error estimate for the lowest-order characteristics-mixed FEM has not been presented although numerous effort has been made in last several decades. In all previous works,  only first-order accuracy in spatial direction was proved under certain time-step and mesh size restrictions. The main purpose of this paper is to establish optimal error estimates, $i.e.$, the second-order in $L^2$-norm for the concentration and the first-order for the pressure/velocity, while the concentration is more important physical component for the underlying model. For this purpose, an elliptic quasi-projection is introduced in our analysis to clean up the pollution of the numerical velocity through the nonlinear dispersion-diffusion tensor and the concentration-dependent viscosity.  Moreover, the numerical pressure/velocity of the second-order accuracy can be obtained by re-solving the (elliptic) pressure equation at a given time level with a higher-order approximation. Numerical results are presented to confirm our theoretical analysis. 

\end{abstract}


\section{Introduction}
\setcounter{equation}{0}
In many engineering areas, one often solves the following miscible displacement system modeling 
an incompressible flow in a porous medium
$\Omega$ 
\begin{align}
&
\Phi\frac{\partial c}{\partial t}-\nabla\cdot(D(\mathbf{u})\nabla c)
+\mathbf{u}\cdot\nabla c=c_1 q^I-c q^P, 
\label{e1}\\
&
\nabla\cdot\mathbf{u}=q^I-q^P, 
\label{e2}\\
&
\mathbf{u}=-\frac{k(x)}{\mu(c)}\nabla p, 
\label{e3}
\end{align} 
for $t\in[0,T]$, with the initial condition
\begin{align}
c(x,0)=c_0(x), \qquad 
\textrm{for }x\in\Omega, \label{e4}
\end{align}
where we assume that the domain $\Omega \in {\bf R}^d$, $d=2,3$,  is
bounded and the condition $\int_{\Omega}pdx=0$ is enforced for the uniqueness
of the solution. In the above system, $c$ represents the concentration of one of the
fluids, $\mathbf{u}$ the Darcy velocity and $p$ the pressure of the fluid mixture.
$\Phi$ denotes the porosity of the medium,
$q^I$ and $q^P$ are given injection and production sources,
$c_1$ is the concentration of the first component in the injection source,
$D(\mathbf{u})=[D_{ij}(\mathbf{u})]_{d\times d}$ is the
diffusion-dispersion tensor (see \cite{BB} for details),
$k(x)$ is the permeability of the medium and
$\mu(c)$ is the concentration-dependent viscosity of the fluid mixture.

In the last several decades, numerical methods and analyses for the
miscible displacement system \refe{e1}-\refe{e3} have been studied 
extensively, $e.g.$, see 
 \cite{DEPT, Ewing, LWS, Pea} and references therein.
Two review articles were written by Ewing and Wang \cite{EWang} and Scovazzi et. al. \cite{SWML}, 
respectively.  In particular, 
Ewing and Wheeler \cite{EW} proposed a fully discrete Galerkin-Galerkin finite
element method for the miscible displacement problem in two dimensional space.
Later, Douglas et al. \cite{DEW} introduced a Galerkin-mixed
finite element method for solving the system \refe{e1}-\refe{e3}.
In both \cite{DEW} and \cite{EW}, a linearized semi-implicit Euler scheme
was applied for the time discretization, and a time step condition
$\tau=o(h)$ was required to obtain optimal error estimates.
Since the concentration equation \refe{e1} is often convection-dominated,
$i.e.$, the diffusion coefficient $D$ is small in many applications,
the characteristics time discretization is more effective for solving this system.
A modified method of characteristics (MMOC) with both finite difference and
finite element approximations was proposed by Douglas and Russell
\cite{DR} for linear convection-dominated diffusion problems. The method is based on the backward Euler scheme in the characteristic time direction and 
classical Galerkin FE approximations in spatial direction. 
The method was extended to the nonlinear miscible displacement equations
in \cite{ERW} with a Galerkin-mixed approximation, where the error estimate 
\begin{align}
  \|c^n-c_h^n\|_{L^2}+\|p^n-p_h^n\|_{L^2}+\|\u^n - \u_h^n\|_{H(div)}
  \le C (\tau + h_c^{r+1} + h_p^{k+1})
  \,
  \label{gm-error}
\end{align}
was established for $d=2$ under the time step restriction $\tau = o(h_p)$ and some mesh size conditions, where 
$D(\u)$ is assumed to be global Lipschitz satisfying 
\begin{align} 
\frac{\partial D(x, {\bf v} )}{\partial {\bf v} } \le K^* 
\label{Du-cond} 
\end{align}    
and $h_c$ and $h_p$ denotes the mesh size of the partition for the concentration 
equation and the pressure equation, respectively. 

The most commonly-used Galerkin-mixed method in practical computation 
is the lowest order one (k = 0,r = 1) \cite{CCW,CWW, DEW, Dur, EW, GS, SWML,SW}.
For the lowest-order mixed method, the error estimate 
\refe{gm-error} reduces to
\begin{align} 
  \| c^n - c_h^n \|_{L^2}  +  \|p^n - p_h^n \|_{L^2}
  + \|\u^n - \u_h^n \|_{H(div)} 
  \le C (\tau + h_p + h_c^2)
  \label{error-0} 
\end{align}
under the more tightened restriction 
\begin{align}
  \tau \le O(h_p^2) 
\, , 
  \label{cond-1}
\end{align} 
for $d=3$ (see (4.42) in \cite{ERW}).   
Numerous effort has been devoted to weakening the time step restriction and mesh size condition 
\cite{DRW,Dur, Rus,SY,WSS}. 
Amongst them, Duran \cite{Dur} showed the error estimate \refe{gm-error}
 under a weaker time-step 
restriction $\tau = o(h_c)$ for $d=3$, the Lipschitz condition \refe{Du-cond} 
for $D(\u)$ and $k\ge 1$. Analysis can be extended to the case $k=0$ as pointed out by the author.     
Further improvement was given recently in \cite{WSS}, 
where in terms of an error splitting technique the 
above error estimate was proved almost unconditionally, i.e., under 
the condition  $\tau \le o(1)$ 
and without the Lipschitz condition \refe{Du-cond} for $D(\u)$. However, the 
analysis was limited to $k =r \ge 1$ which exclude the popular lowest-order mixed method. 
Moreover the modified method of characteristics combined with 
many other approximations in spatial direction has also been studied extensively  \cite{DRW,DO,FN,KHR,KY,LWC,MLY,SY}. 
To maintain the conservation of the mass, a related Eulerian-Lagrangian localized adjoint
method (ELLAM) was studied in \cite{CRHE, Wang} for advective-diffusive equations, in which 
an ELLAM scheme was used in time direction. 
Analysis of an ELLAM-MFEM for \refe{e1}-\refe{e3} was presented in \cite{Wang,WER}. 
 A more general ELLAM scheme was proposed and investigated in the recent work \cite{CDL}. 
The convergence rate of the method in spatial direction is similar to those in \refe{gm-error} and 
\refe{error-0}. 
Some other type methods of characteristics can be found in
\cite{ASW,KM,MP}. 
In addition, the characteristics type methods have been applied and analyzed for
many other linear and nonlinear parabolic PDEs from various engineering applications
\cite{ASW,BNV,DEPT,FN,GPP,LWC,SWS}. Numerical simulations show that
the time-truncation errors of the MMOC are much smaller
than those of standard schemes for convection-dominated models.

There are still several issues to be further addressed 
for the popular lowest-order characteristics-mixed FEM.
(i). In the lowest-order characteristics-mixed method, 
a linear Lagrange approximation and a zero-order Raviart-Thomas approximation 
are used for the concentration and the pressure/velocity, respectively. Clearly,  
the error estimate presented in \refe{error-0} is not optimal for the concentration in $L^2$-norm, 
while the concentration is a more important physical component in practical applications.  
(ii). The modified method of characteristics is based on characteristic tracking, along which 
the method may greatly reduce the temporal error
and allow one to use a large time step in computations.
However, certain tightened time-step condition was always required in previous analysis. 
(iii). The Lipschitz condition \refe{Du-cond}  for the diffusion-disperson tensor $D(\u)$ 
may not be realistic in practice.  Analysis of Galerkin-mixed 
FEMs for \refe{e1}-\refe{e3} under a weaker assumption of $D(\u)$ being smooth 
(without the global Lipschitz condition \refe{Du-cond}) 
was done in \cite{CCW,DEW}, which, however, leads to some more serious mesh condition 
for the lowest-order mixed method.

This paper focuses on a new analysis 
of the lowest-order characteristics-mixed finite element
method for the nonlinear and coupled system \refe{e1}-\refe{e3}. 
We shall establish the optimal  $L^2$-norm error estimates 
\begin{align}
 & \| c^n - c_h^n \|_{L^2}  \le C (\tau + h^2) 
 \label{error-1} 
 \\ 
 &  \|p^n - p_h^n \|_{L^2}
  + \|\u^n - \u_h^n \|_{H(div)}  
  \le C (\tau + h)
  \label{error-2}
\end{align}
only under the condition 
\begin{align} 
\tau = o\left ( \frac{1}{|\log h|} \right ) 
\label{t-cond}
\end{align} 
and the weak assumption of 
$D(\u)$ being smooth (without the global Lipschitz condition \refe{Du-cond}).  
The new analysis shows 
that the method provides a second-order accuracy for the concentration, while 
only first-order accuracy was proved in previous works. Moreover, 
with the numerical concentration of second-order accuracy, a second-order 
pressure/velocity at a given time level can be obtained by re-solving the elliptic pressure equation 
with a first-order RT approximation. The extension to more general cases with higher-order 
approximations and different mesh partitions can be made analogously. 
The analysis presented in this paper 
is based on an elliptic quasi-projection proposed in \cite{SW}, 
an error splitting technique presented in \cite{LS2} and negative-norm estimate of the numerical 
velocity. 
With the quasi-projection and some more precise estimates in the characteristic direction, 
the lower-order approximation to the pressure/velocity does not pollute the accuracy of numerical 
concentration in our analysis. 
The optimal analysis under the weaker time step condition \refe{t-cond} is given in terms of the error splitting technique.

The paper is organized as follows. In Section 2, we present 
our notations and our main results. A new re-covering technique is introduced, with which 
the second-order accuracy of numerical velocity/pressure can be obtained by re-solving the 
elliptic pressure equation with a higher-order approximation and the obtained numerical concentration. In section 3, we first present several useful lemmas and more precise estimates in the 
characteristic direction. In terms of  the error splitting argument, we analyze the temporal and 
spatial errors, respectively and 
the boundedness of numerical solutions. 
Then, we present optimal error estimates of the
numerical scheme. In Section 4,
numerical results are given to confirm our theoretical analysis.

\section{Main results}
\setcounter{equation}{0}
We at first define some notations used in this paper.
For any integer $m\geq 0$ and $1\leq p\leq\infty$, let
$W^{m,p}(\Omega)$ be the Sobolev space of functions with the norm
\begin{align*}
\|f\|_{W^{m,p}}=\left\{
\begin{array}{ll}
\Big(\displaystyle\sum\limits_{|\beta|\leq m}
\int_{\Omega}^{}|D^{\beta}f|^pdx\Big)^{\frac{1}{p}},
&\textrm{for }1\leq p<\infty,\\
\displaystyle\sum\limits_{|\beta|\leq m}\textrm{ess}\sup_{\Omega}|D^{\beta}f|,
&\textrm{for }p=\infty,
\end{array}\right.
\end{align*} 
where
\begin{eqnarray*}
D^{\beta}=\frac{\partial^{|\beta|}}{\partial x_1^{\beta_1}\cdots\partial x_d^{\beta_d}},
\end{eqnarray*} 
for the multi-index $\beta=(\beta_1,\cdots,\beta_d)$,
$\beta_1\geq 0,\cdots,\beta_d\geq 0$, and $|\beta|=\beta_1+\cdots+\beta_d$.
When $p=2$, we denote $W^{m,2}(\Omega)$ by $H^m(\Omega)$.
We define
$L_0^k(\Omega)=\{f\in L^k(\Omega):\int_{\Omega}f dx=0\}$ and
$H(\textrm{div};\Omega)=\{\mathbf{f}=(f_1,\cdots,f_d):f_i,
\nabla\cdot\mathbf{f}\in L^2(\Omega),1\leq i\leq d\}$. 
For simiplicity, we write ${\bf f} \in W^{m,p}$ if $f_i \in W^{m,p}$. 
To avoid technical difficulties on boundary, we assume
that $\Omega$ is a rectangle in $\mathbb{R}^2$
(or cuboid in $\mathbb{R}^3$) and the problem \refe{e1}-\refe{e3}
and the corresponding FE spaces are $\Omega$-periodic as usual 
\cite{ERW,MLY,Rus,WSS}.

Let $\pi_h$ be a quasi-uniform partition of $\Omega$ into
triangles $T_j$, $j=1,\cdots, M$, in $\mathbb{R}^2$ (or
tetrahedra in $\mathbb{R}^3$) of diameter less than $h$.
We denote  by $(S_h^k, H_h^k)$ $k$-order Raviart-Thomas finite element space \cite{RT} 
\begin{align} 
& S_h^k : = \{ w \in L^2_0(\Omega): w|_{T_j} \in P_{k} \} 
\nn \\ 
& H_h^k : = \{ \v \in H(\textrm{div};\Omega): \v|_{T_j} \in P_k \otimes xP_k \}  
\nn 
 \end{align} 
 and 
 by $V_h^1$ the standard linear Lagrange FE space on the partition $\pi_h$
 where $P_k$ denotes the polynomial space of degree $\le k$. 
 
Let $\{t_n|t_n=n\tau;0\leq n\leq N\}$ be a uniform
partition of $[0,T]$ with the time step $\tau=T/N$, and we denote
\begin{eqnarray*}
c^n(x)=c(x,t_n),\quad
\mathbf{u}^n(x)=\mathbf{u}(x,t_n),\quad
p^n(x)=p(x,t_n).
\end{eqnarray*} 
For a sequence of functions
$\{\omega^n\}_{n=0}^N$, we define
\begin{eqnarray*}
D_{\tau}\omega^{n+1}=\frac{\omega^{n+1}-\omega^n}{\tau}.
\end{eqnarray*}

Here we assume that
the permeability $k(\cdot)$ is in the space $H^2(\Omega)$ 
satisfying 
\begin{align} 
& k_0^{-1}\leq k(x)\leq k_0\quad \textrm{for }x\in\Omega 
\nn 
\end{align} 
and 
the concentration-dependent viscosity
$\mu(\cdot) \in H^2(\mathbb{R})$ is globally Lipschitz, satisfying 
\begin{align} 
& \mu_0^{-1}\leq \mu(x)\leq \mu_0\quad \textrm{for }x\in\Omega 
\label{mu} 
\end{align} 
for some positive constants $k_0$ and $\mu_0$.
Moreover, the injection and production sources satisfy
\begin{eqnarray}
\|q^I\|_{W^{1,4}},\|q^P\|_{W^{1,4}}\leq K_1.
\end{eqnarray}
The diffusion-dispersion tensor
$D(\u)=\Phi(d_{mt}(|\u|) I+d_{lt}(|\u|) \u \otimes \u)$ is a $d\times d$ matrix,
where $d_{mt}(z) > d_m>0$, $d_{lt}(z) > 0$ for $z>0$ and $\u \otimes \u = \u \u^T$. We further assume that $d_{mt}(z), d_{lt}(z) \in H^3(R)$. But $D(\u)$ may not be globally Lipschitz. 
 For the system \refe{e1}-\refe{e3} being well-posed,
we add
\begin{align}
\int_{\Omega}q^I dx=\int_{\Omega}q^Pdx.
\end{align}
For simplicity, we assume that $\Phi = 1$. 
These assumptions have been made in those previous analysis as usual 
\cite{Dur, ERW, EW, LS2, LS3,WSS}.

With the above notations, the modified method of characteristics with
the lowest-order mixed FE approximation is to find $(c_h^n, p_h^n, \u_h^n ) \in (V^1_h, S_h^0, H_h^0)$
such that
\begin{align} 
& \left( \frac{c_h^{n+1}-c_h^n(x_{\u_h^n})}{\tau},\phi_h \right)
+\Big(D(\u_h^n)\nabla c_h^{n+1},\nabla\phi_h\Big)
=\Big(c_1q^I-c_h^{n+1}q^P,\phi_h \Big), 
\label{c1}
\\ &
\left( \frac{\mu(c_h^n))}{k(x)}  \u_h^{n}, \v_h \right)
=\Big(p_h^{n},\nabla\cdot \v_h\Big), 
\label{c2}  
\\ &
\Big(\nabla\cdot \u_h^{n},\varphi_h \Big)
=\Big(q^I-q^P,\varphi_h\Big),  
\label{c3} 
\end{align} 
for all $(\phi_h, \varphi_h, \v_h) \in (V^1_h, S_h^0, H_h^0)$, 
where
\begin{align*} 
x_{\u_h^n}(x): =x-\u_h^n(x)\tau, \qquad \mbox{ for $x \in \Omega$} 
\end{align*} 
and $c_h^0=I_h c_0$ with $I_h$
being the Lagrangian interpolation operator. 
Some slightly different schemes were investigated by many authors 
 \cite{DRW,Dur,ERW,Rus,WSS}. Error estimates of all these schemes were obtained with some restrictions on time step and spatial mesh size and under certain assumptions for the diffusion-dispersion tensor $D(\u)$. 
 It is easy to extend our analysis to these schemes. 

For simplicity, here we assume that the system \refe{e1}-\refe{e3}
admits a unique solution satisfying
\begin{align}
\|c_0\|_{H^{2}}
+\|c\|_{L^{\infty}(I;H^{2})}
+\|c_t\|_{L^{\infty}(I;H^{2})}
+\|c\|_{L^{\infty}(I;W^{2,4})}
+\|c_{tt}\|_{L^{2}(I;L^2)}
&\nn\\
+\|\mathbf{u}\|_{L^{\infty}(I;W^{2,4})}
+\|\mathbf{u}_t\|_{L^{\infty}(I;L^2)}
+\|\mathbf{u}_t\|_{L^{2}(I;H^1)}
+\|p\|_{L^{\infty}(I;H^2)}
&\leq K_2.  
\label{reg}
\end{align}
Theoretical analysis for the underlying system can be found in 
\cite{Feng}. 
The present paper focuses on the optimal error estimates of the lowest-order characteristics-mixed FEM, while the above regularity assumptions may be weakened slightly. 

Next we present our main results in the following theorem.
\begin{theorem}\label{main}
{\it Suppose that the system \refe{e1}-\refe{e4}
has a unique solution $(c,\mathbf{u},p)$ satisfying \refe{reg}. 
Then, there exists a positive constant $h_0$  
such that when $h<h_0$, the finite element system \refe{c1}-\refe{c3} admits a unique solution $(c_h^m, \u_h^m, p_h^m) \in (V^1_h, S_h^0, H_h^0)$,
$m=0,1,\cdots,N$. Moreover, under the condition 
\begin{align} 
\tau = o\left ( \frac{1}{|log h| } \right ) , 
 \label{t-cond} 
 \end{align}  
the FE solution satisfies 
\begin{align} 
& \max_{0\leq m\leq N} \|c_h^m-c^m\|_{L^2} 
\le C_0 (\tau + h^2) 
\label{error0} 
\\
& \max_{0 \le m \le N} \left ( \|\u_h^m-\u^m\|_{H(\textrm{div})}
+ \|p_h^m-p^m\|_{L^2} \right ) 
\leq C_0(\tau+h) 
\label{error}
\end{align} 
where $C_0$ is a positive constant independent of $m$, $\tau$ and $h$ 
and may be dependent upon $K_2$ and the physical constants,
$K_1$, $k_0$ and $\mu_0$.
}
\end{theorem}

With the obtained 
numerical solution $(c_h^n, \u_h^n, p_h^n)\in (V^1_h, S_h^0, H_h^0)$, a new numerical velocity/pressure 
of a second-order accuracy  
can be obtained by re-solving the pressure equation  
\begin{align}
& \left( 
 \frac{\mu(c_h^n)}{k(x)}  \widehat \u_h^{n}, \v_h \right)
=\Big(\widehat p_h^n, \nabla\cdot \v_h\Big),\qquad 
\v_h\in H_h^1 
\label{new-p} 
\\ &
\Big(\nabla\cdot \widehat \u_h^{n},\varphi_h \Big)
=\Big(q^I-q^P,\varphi_h\Big)
, \qquad \varphi_h\in S_h^1,
\label{new-u} 
\end{align} 
with the first-order mixed FE approximation $(\widehat p_h^n, \widehat \u_h^n) \in (S_h^1, H_h^1)$ at a given time level $t_n$. 

\begin{corollary}\label{u-Error}
{\it Suppose that the system \refe{e1}-\refe{e4}
has a unique solution $(c,\mathbf{u},p)$ satisfying \refe{reg} . The FE 
solution $(\widehat p_h^n, \widehat \u_h^n) \in (S_h^1, H_h^1)$ of the system 
\refe{new-p}-\refe{new-u} 
satisifies 
\begin{align}
&\|\widehat \u^n_h - \u^n\|_{L^2} + \|\widehat p_h^n - p^n \|_{L^2}  \le \widehat C_0 (\tau+h^2),   
\label{error-u}
\end{align}
where  $\widehat C_0$ is a constant independent of $n$, $h$ and $\tau$ 
and may be dependent upon $K_2$, $C_0$ and the physical constants,
$K_1$, $k_0$ and $\mu_0$.
}
\end{corollary}

In the rest of this paper, we denote by $C$ a generic positive constant and by $\epsilon$ a 
generic small positive constant, which are independent of $n, h,\tau, C_0$ and $\widehat C_0$.  
The following classical 
Gagliardo-Nirenberg inequality \cite{N} will be frequently used in our proof, 
\begin{align}
\|\partial^j u\|_{L^p}\leq
C\|\partial^m u\|_{L^r}^{\alpha}\|u\|_{L^q}^{1-\alpha} + C\|u\|_{L^q},
\label{GN} 
\end{align} 
for $0\leq j<m$ and $\frac{j}{m}\leq\alpha\leq 1$
with
$$ 
 \frac{1}{p}=\frac{j}{d}+\alpha\left(\frac{1}{r}-\frac{m}{d}\right)
+(1-\alpha)\frac{1}{q}, 
$$  
except $1<r<\infty$ and $m-j-\frac{ d}{r}$ is a
non-negative integer, in which case the above estimate holds only
for $\frac{j}{m}\leq\alpha<1$. Moreover, we present a classical 
    discrete Gronwall's inequality in the following lemma. 
    \begin{lemma}\label{l2-1}
{\it  
    Let $\tau$, $B$ and $a_{k}$, $b_{k}$, $c_{k}$, $\gamma_{k}$, 
    for integers $k \geq 0$, be non-negative numbers such that
    \[
      a_{n} + \tau \sum_{k=0}^{n} b_{k} 
      \leq \tau \sum_{k=0}^{n} \gamma_{k} a_{k} + 
      \tau \sum_{k=0}^{n} c_{k} + B \, , \quad \mathrm{for } \quad n \geq 0 \, ,
    \]
    suppose that $\tau \gamma_{k} < 1$, for all $k$,
    and set $\sigma_{k}=(1-\tau \gamma_{k})^{-1}$. 
    Then
    \[
      a_{n} + \tau \sum_{k=0}^{n} b_{k} 
      \leq  \exp(\tau \sum_{k=0}^{n} \gamma_{k} \sigma_{k}) (\tau \sum_{k=0}^{n} c_{k} + B) \, , 
      \quad \mathrm{for } \quad n \geq 0 \, .
    \]
} 
\end{lemma}

\section{Analysis}   
\setcounter{equation}{0}
Before proving our main theorem, we present several lemmas in the following subsection, which are useful in the proof of the main theorem.  
\subsection{Prelimaries} 
\begin{lemma}\label{l3-1}
{\it 
Assume that $f\in L^p(\Omega)$ is 
$\Omega$-periodic and 
$g$  is a piecewise smooth function satisfying 
\begin{align}
\tau | g(x_a) - g(x_b) | \le  \frac{1}{2} ( |x_a - x_b| + h ), \qquad \mbox{ for any } x_a, x_b \in \Omega . 
\label{Lip}  
\end{align} 
Then 
\begin{align} 
\| f(x+\tau g(x)) \|_{L^p} \le C \| f(x) \|_{L^p} \, . 
\label{Lip-int} 
\end{align} 
} 
\end{lemma} 

{\it Proof.} A special case of \ref{Lip-int} was studied in \cite{ERW}. 
Letting $z_g(x) = x + \tau g(x)$, 
by \refe{Lip} we have 
$$ 
| z_g(x_a) - z_g(x_b) | = | (x_a -x_b) + \tau ( g(x_a) - g(x_b)) | 
\ge \frac{1}{2} (|x_a - x_b| - h) 
$$ 
which shows that 
$z_g(x_a)$ and $z_g(x_b)$ are not in one element when 
$|x_a - x_b| >3h$. Hence $z_g(x)$ is  globally at most finitely-many-to-one 
and maps $\Omega$ into itself and its immediate-neighbor periodic copies. 
By noting 
$$ 
\| f(x+\tau g(x) \|_{L^p}^p = \sum_{j=1}^M \int_{T_j} |f(z_g)|^p dx \, ,  
$$ 
we see that the sum above is bounded by 
finitely many multiples of the integral $\int_{\Omega} |f(x)|^p dx$ \cite{ERW}. 
\refe{Lip-int} follows immediately. 
\quad \endproof
\vskip0.1in 

Clearly, \refe{Lip} holds 
if $g \in W^{1,\infty}(\Omega)$ or $\tau \| g \|_{L^\infty} \le h/4$. 
Analysis for the method of characteristics type relies on the approximation in the 
characteristic direction. Several estimates along the characteristic direction were presented in \cite{Dur,ERW,WSS}. 
In the following lemma we present some more precise estimates, which play an important role in our analysis.

\begin{lemma}\label{l3-2}
{\it  
Assume that  $v, \rho \in C^0(\Omega)\cap H^1(\Omega)$,  
$g_1, g_2$  are $\Omega$-periodic and piecewise smooth and 
$g_2, (g_1-g_2)$ 
satisfy the condition \refe{Lip}. 
Then (i)  we have 
\begin{align}
| \left ( \rho(x-g_1\tau)-\rho(x-g_2\tau), \, v \right )|  
\leq
C\tau\|\rho\|_{W^{1,p}} \|g_1-g_2 \|_{L^q} \| v \|_{L^6} 
\label{l2-2-1} 
\end{align}
where $1/p+1/q=5/6$;  
 (ii) 
 if  $g_1, g_2 \in W^{2,3}(\Omega)\cap C^1(\Omega)$,  
\begin{align}
| \left (  \rho(x-g_1\tau)-\rho(x-g_2\tau), v \right ) | 
\leq
C\tau\| \rho\|_{L^p}  \|g_1-g_2 \|_{W^{1,q}} \| v \|_{H^1} 
\label{l2-2-2} 
\end{align}
where $1/p+1/q=1/2$ for $2 \le p < 6$ and (iii) if $\rho \in W^{2,4}(\Omega)$, 
we have 
\begin{align} 
|( \rho(x-g_1\tau)-\rho(x-g_2\tau), \, v)) |  
\le  \tau C 
\| \rho \|_{W^{2,4}} (\| g_1-g_2 \|_{H^{-1}}  + \tau \| g_1 -g_2 \|_{L^4}^2 ) 
\| v \|_{H^1} \, .
\label{l2-2-3}  
\end{align} 
}
\end{lemma}

{\it Proof.}  (i). It is easy to see that 
\begin{align} 
 | (\rho(x-g_1\tau) - \rho(x-g_2\tau)), \, v) | 
 &  =
\left | \left ( \int_0^1 \partial_s \rho(x-g_2\tau 
- s\tau(g_1 -g_2) {\color{red} )}ds, v\right ) \right | 
\nn \\
& = \tau \left | \int_0^1 \int_{\Omega}   \nabla \rho(z(x)) 
{\color{red} \cdot }(g_1-g_2) v(x)  dx ds \right | 
 \nn \\
 & \le \tau \int_0^1  \|\nabla \rho(z(x)) \|_{L^p}  \|  g_1-g_2 \|_{L^q} \|  v \|_{L^6}  ds  
 \label{rho1}
 \end{align} 
 where $z(x) = x - \tau g_2 - s \tau (g_1-g_2)$ defines a map. 
 
Since $\nabla z = I - \tau \nabla g_2 - s\tau (\nabla g_1 - \nabla g_2)$ 
and $g_2, g_1-g_2$ satisfy \refe{Lip}, for any $x \in T_j$, 
$\mbox{det}(\nabla z) >1/2$ 
and 
$$ 
 \| \nabla \rho(z(x)) \|_{L^p}  
  \le C  \|\nabla \rho(x) \|_{L^p}  \, . 
 $$  
 \refe{l2-2-1} follows immediately. 
 
 (ii). For $\rho, v \in C^0(\Omega)\cap H^1(\Omega)$ and $g_1, g_2 \in W^{2,3}(\Omega)\cap C^1(\Omega)$,  we have 
\begin{align} 
\nabla_x \cdot \left [ 
v(g_1-g_2) \cdot (\nabla z)^{-1} \rho(z) \right ] 
&  = v(g_1-g_2) \cdot (\nabla z)^{-1} \cdot  \nabla_x \rho(z(x)) 
\nn \\ 
& + \nabla_x \cdot [v(g_1-g_2) \cdot (\nabla z)^{-1} ]   \rho(z) 
\nn \\ 
& = v(g_1-g_2) \cdot \nabla \rho(z) + \nabla_x  \cdot [v(g_1-g_2) \cdot (\nabla z)^{-1} ]   \rho(z) 
\, . 
\end{align} 
Therefore, by \refe{rho1} 
\begin{align}
 | (\rho(x-g_1\tau) - \rho(x-g_2\tau)), \, v) | 
 & = \tau \left | \int_0^1 \int_{\Omega}   \nabla \rho(z)
 {\color{red} \cdot } (g_1-g_2) v(x)  dx ds \right |  
 \nn \\ 
 & = \tau 
  \left | \int_0^1  \int_{\Omega}  \nabla_x \cdot [v(g_1-g_2) ) \cdot (\nabla z)^{-1}] \rho(z)  dx ds \right | 
 \nn \\ 
&\le C\tau 
  \int_0^1 \|  \rho(z(x)) \|_{L^p(\Omega)}
   \| \nabla_x \cdot [v(g_1-g_2) \cdot (\nabla z)^{-1} ] 
   \|_{L^{q_1} (\Omega)} \, ds 
 \end{align} 
where $1/p+1/q_1=1$. 
Since $g_1, g_2 \in C^1(\Omega)$ and 
 $\rho \in C^0(\Omega)$ are
  $\Omega$-periodic, by Lemma \ref{l3-1} we have 
 \begin{align} 
 \|  \rho(z(x)) \|_{L^p(\Omega)} \le C  \|  \rho(x) \|_{L^p(\Omega)}
 \nn 
 \end{align} 
and by Gagliardo-Nirenberg inequality, 
 \begin{align} 
 \| \nabla_x \cdot [v(g_1-g_2) \cdot (\nabla z)^{-1} ] 
 \|_{L^{q_1}(\Omega)}
& \le 
  \|  v \|_{H^1}  \| g_1-g_2 \|_{L^\infty}  \| z  \|_{W^{1,q}}  
 \nn \\ 
 & +   \|  v \|_{L^6} \| g_1-g_2 \|_{W^{1,q}}  \| z  \|_{W^{1,3}} 
  +  \|  v \|_{L^6} \| g_1-g_2 \|_{L^\infty}  \| z  \|_{W^{2,3}} 
  \nn \\ 
  & \le C \|  v \|_{H^1} \| g_1-g_2 \|_{W^{1,q}} 
  \nn 
  \end{align} 
where we have noted $1/2+1/q = 1/q_1$, $q>3$ and $z(x) \in W^{2,3}(\Omega)\cap C^1(\Omega)$. 
It follows that 
$$ 
| (\rho(x-g_1\tau) - \rho(x-g_2\tau)), \, v) | 
\le C \tau  \| \rho(x) \|_{L^p}  \| v \|_{H^1} \| g_1 -g_2 \|_{W^{1,q}} \, .
$$  
We have proved \refe{l2-2-2}. 
 
(iii). 
Since 
\begin{align}
 \rho(x-g_1\tau) - \rho(x-g_2\tau) & = -\tau \nabla \rho(x-\tau g_2) \cdot (g_1 -g_2) 
\nn \\ 
&  + \frac{1}{2} \tau^2 (g_1-g_2) \cdot \int_0^{\bar s} \nabla^2 \rho(x - \tau g_2 - s\tau (g_1-g_2)) 
 \cdot  (g_1-g_2) ds
\nn 
\end{align} 
for some $\bar s$ with $0 < \bar s <1$, we get    
\begin{align} 
 | (\rho(x-g_1\tau) - \rho(x-g_2\tau)), \, v) | 
 & \le  \tau \|  \nabla \rho(x-\tau g_2) v(x)\|_{H^1} 
  \| g_1-g_2 \|_{H^{-1}}  
\nn \\ 
& +  \tau^2  \int_0^1 \| g_1- g_2 \|^2_{L^4} \| v \|_{L^6} 
 \| \nabla^2 \rho(z(x)) \|_{L^3} ds 
 \label{rho2}
 \end{align} 
Similarly we have 
\begin{align} 
&  \| \nabla^2 \rho(z(x)) \|_{L^3} 
 \le C\|  \rho(x) \|_{W^{2,3}} 
\nn \\ 
& \|  \nabla \rho(x-\tau g_2) v(x)\|_{H^1}  \le C
 \|  \nabla \rho \|_{W^{1,4}} \| v \|_{H^1} \, . 
 \nn 
\end{align} 
\refe{l2-2-3} follows immediately. 
 The proof is complete. \quad \endproof 
 \vskip0.1in

To prove Theorem 2.1, we introduce a characteristic time-discrete system:
\begin{align} 
& \frac{\C^{n+1}-\C^n(x_{U^n})}{\tau}
-\nabla\cdot(D(U^{n})\nabla \mathcal{C}^{n+1}) 
=c_1q^I-\mathcal{C}^{n+1}q^P, 
\label{t-e1}
\\ 
& U^n=- \frac{k(x)}{\mu(\C^n)} \nabla P^n, 
\label{t-e2}\\
& \nabla\cdot U^n=q^I-q^P, 
\label{t-e3}
\end{align} 
with periodic boundary conditions and the following initial condition
\begin{align} 
{\cal C}^0(x)=c_0(x),
\nn 
\end{align} 
where $x\in\Omega$, $t\in[0,T]$ and
\begin{align*} 
x_{U^n}(x): =x-U^n(x) \tau.
\end{align*}
The condition $\int_{\Omega}P^{n}dx=0$ is enforced
for the uniqueness of the solution.
The above system can be viewed as an iterated sequence of elliptic PDEs and the numerical solution 
$(c_h^n, p_h^n, \u_h^n)$ can be viewed as the FE solution of the elliptic system \refe{t-e1}-\refe{t-e3}. 
We present the regularity of the solution of the system \refe{t-e1}-\refe{t-e3} 
and the corresponding error estimates 
in the following lemma. The proof is omitted since a slightly 
different lemma was proved in \cite{WSS}. 

\begin{lemma}\label{l3-3}
{\it Suppose that the system \refe{e1}-\refe{e4}
has a unique solution $(c,\u,p)$ satisfying \refe{reg}.
Then, there exists $\tau_1>0$ such that
when $\tau<\tau_1$, the time-discrete system \refe{t-e1}-\refe{t-e3}
admits a unique solution $(\mathcal{C}^n, U^n, P^n)$,
$n=0,1,\cdots,N$, which satisfies
\begin{align}
&\|c^n - \C^n \|_{H^1} + \|\u^n - U^n\|_{H^1}
+ \|p^n - P^n \|_{H^1} \leq C_1\tau , 
\label{reg-1}
\\
& \|U^n\|_{W^{2,4}} + \| {\cal C}^n \|_{W^{2,4}} 
+ \sum_{m=1}^n \tau\|D_{\tau}\mathcal{C}^m\|_{H^2}^2
\leq C_1 
\label{reg-2}
\end{align}
where $C_1$ is a constant independent of $h$, $\tau$, $n$, $C_0$
and may depend upon $K_1$, $K_2$, $k_0$ and $\mu_0$.
}
\end{lemma}

Moreover, for any fixed integer $n\geq 0$, we denote
by $(\wt P^{n}_h, \wt  U_h^{n})$ the mixed projection of
$(P^{n}, U^{n})$ on
$S_h^0 \times \H_h^0$ such that
\begin{align}
  & \biggl(\frac{\mu(\C^{n})}{k(x)} (\widetilde{U}_h^{n}-U^{n}),\,\v_h\biggl)
  = \Big({\widetilde  P}_h^{n} - P^{n}  ,\, \nabla \cdot \v_h \Big),
  \label{e-proj-1}
  \\[3pt]
  & \Big(\nabla\cdot ({\wt  U}_h^{n}-U^{n}) ,\, \varphi_h\Big) =0, \quad 
  \forall (\varphi_h,\v_h)
  \in S_h^0 \times \H_h^0 
  \label{e-proj-2}
\, . 
\end{align}
Error estimates of the mixed projection are presented below. 
\begin{align}
  &\|U^n-\widetilde  U_h^n\|_{L^p} + \|P^n- \widetilde  P_h^n\|_{L^p} 
  +  \| U^n - \widetilde U_h^n \|_{H(\textrm{div})}     
   \le C h,   \quad \mbox{for~all}~~2 \le p \le 4,   
  \label{Up-Lp-div} \\
  & \|U^n-\widetilde  U_h^n\|_{L^\infty} \leq Ch \log(1/h) 
  \label{U-infty} \\
  &\|U^n- \widetilde  U_h^n\|_{H^{-1}} + \|P^n- \widetilde  P_h^n\|_{H^{-1}}  \leq Ch^2 \, . 
  \label{Up-H-1} 
\end{align}
The proof of \refe{Up-Lp-div} follows classical mixed FE theory \cite{BS, Dur2, RT} and the proof of 
\refe{U-infty}-\refe{Up-H-1} can be found in 
in \cite{Dur2,GN} and \cite{DRob}, respectively.

For a given $U^n$, an elliptic quasi-projection $\wt \C_h^{n+1}$ of $\C^{n+1}$ from $H^1(\Omega) \rightarrow
V_h^1$ is defined by
\begin{align}
\label{Nonclassical-0}
  \Big(D( \wt U^n_h)\nabla \wt \C^{n+1}), \, \nabla \phi_h \Big) 
=  \Big(D( U^n)\nabla \C^{n+1}), \, \nabla \phi_h \Big),   
  \qquad \mbox{for~all}~~\phi_h\in V_h,~~n\geq 0, 
\end{align}
with $\int_\Omega(\wt \C^{n+1}_h-\C^{n+1})d x=0$ and $\wt \C^0_h = I_h\C^0$.  
The above equation is equivalent to 
\begin{align} 
\label{Nonclassical} 
  \Big(D( U^n)\nabla ( \wt \C^{n+1}_h - \C^{n+1}), \, \nabla \phi_h \Big) 
  + \Big((D(\wt U_h^n)-D(U^n))\nabla \widetilde \C^{n+1}_h ,  
\, \nabla \phi_h \Big)= 0. 
\end{align} 
All previous analyses were based on a classical elliptic projection 
proposed  in \cite{Whe}, where the second term in \refe{Nonclassical} is excluded. Applying the classical elliptic projection for the present nonlinear and strongly coupled problem 
leads to serious pollution in estimating the error of concentration. 
Here the quasi-projection is used in our analysis. 
Some basic estimates of the quasi-projection
are presented in the following lemma and the proof can be found in \cite{SW} 
\begin{lemma}\label{l3-4}
  {\it Under the assumptions of Theorem \ref{main},
    there exists $h_1>0$ such that for any $h\leq h_1$ and $2\le p \le 4$
    \begin{align}
      &\|\C^n- \wt \C^n_h\|_{L^2} 
      +h\|\nabla (\C^n- \wt \C^n_h)\|_{L^p}
      \leq  C_2h^2,   
      \label{p2-1}  
    \end{align}
and
    \begin{align}
      &\biggl(\sum_{n=0}^{N-1}\tau\|D_t (\C^n
      - \wt \C^n_h)\|_{L^2}^2\biggl)^{1/2}
      \leq  C_2h^2
      \label{p2-2}
\end{align} 
where $C_2$ is a constant independent of $h$, $\tau$, $n$, $C_1$
and may be dependent upon
    $K_1$, $K_2$, $C_0$, $k_0$ and $\mu_0$.
}
\end{lemma}

Under the regularity assumption (\ref{reg}), we can see  from Lemma \ref{l3-4} that 
\begin{align} 
\| \wt \C_h^n \|_{W^{1,\infty}} \le \| \C^n \|_{W^{1,\infty}} + \| \wt \C_h^n - \C^n \|_{W^{1,\infty}} 
\le C + C h^{-3/4} \| \wt \C_h^n - \C^n \|_{W^{1,4}}   \le C 
\label{chh} 
\end{align} 
for $n=1,2,...,N$. 

\subsection{The proof of Theorem \ref{main} } 
Since at each time step, 
 the coefficient matrix of the system \refe{c1} is symmetric positive definite 
 and \refe{c2}-\refe{c3} defines  a standard saddle point system,   
 the existence and uniqueness of 
the numerical solution $(c_h^{n+1}, \u_h^n, p_h^{n})$ 
follows immediately. 

The key to the proof of 
\refe{error0}-\refe{error} is the boundedness of numerical solution. 
In terms of temporal-spatial error splitting argument introduced in \cite{LS1},  we have 
\begin{align}
&
\|c_h^n-c^n\|_{L^2}\leq
\|c_h^n - \C^n\|_{L^2} 
+ \| \C^n - c^n \|_{L^2} , 
\nn \\
&
\|\u_h^n-\u^n\|_{L^2}\leq
\|\u_h^n - U^n\|_{L^2}
+ \| U^n - \u^n \|_{L^2},  
\label{splitting} 
\\ 
&
\|p_h^n-p^n\|_{L^2}\leq
 \|p_h^n - P^n\|_{L^2} 
+ \| P^n - p^n \|_{L^2}\, . 
\nn 
\end{align}

By noting Lemma \ref{l3-3}, we only need to estimate the first  
terms in the splitting above. 
To estimate them, we make a further splitting in terms of the mixed  
projection and quasi-projection introduced above to get  
\begin{align} 
& \| c_h^n - \C^n \|_{L^2} \le \| \xi_c^n \|_{L^2} + \| e_c^n  \|_{L^2} 
\nn \\ 
& \| \u_h^n - U^n \|_{L^2} \le \| \xi_u^n \|_{L^2} + \| e_u^n \|_{L^2} 
\label{splitting-2} 
\\ 
& \| p_h^n - P^n \|_{L^2} \le \| \xi_p^n \|_{L^2} + \| e_p^n  \|_{L^2} 
\nn 
\end{align} 
where  
\begin{align*} 
& \xi_c^n=c_h^n- \wt \C^n_h, \quad 
\xi_{u}=\u_h^n- \wt U^n_h,\quad 
\xi_p^n=p_h^n-\wt P^n_h, \quad n=0,1,\cdots,N 
\\ 
&  e_c^n= \C^n - \wt \C^n_h, \quad 
e_u^n= U^n - \wt U^n_h,\quad 
e_p^n=P^n - \wt P^n_h, \quad n=0,1,\cdots,N  \, . 
\end{align*} 
The mixed projection and quasi-projection errors $\e_c^n, e_u^n, e_p^n$ have been presented in \refe{Up-Lp-div}-\refe{Up-H-1} 
and \refe{p2-1}-\refe{p2-2}, respectively. 
Since  Lemma \ref{l3-3} and Lemma \ref{l3-4} have been proved,  hereafter  
we assume that the generic constant $C$ 
may depend upon $C_1$ and $C_2$.

The weak formulation of the characteristic time-discrete system \refe{t-e1}-\refe{t-e3} can be 
written by 
\begin{align} 
& \left (  \frac{\C^{n+1}-\C^n(x_{U^n})}{\tau}, \phi \right ) 
 + \left ( (D(U^{n})\nabla \mathcal{C}^{n+1}),   \nabla \phi \right ) 
= \left ( c_1q^I-\mathcal{C}^{n+1}q^P, \,  \phi \right ) 
\nn \\ 
& \left ( \frac{\mu(\C^n)}{k(x)} U^n, \, \v \right ) = \left ( P^n, \, \nabla \cdot \v \right )  
\nn \\
& ( \nabla\cdot U^n, \psi )= ( q^I-q^P, \psi)  
\nn 
\end{align} 
with periodic boundary conditions and the following initial condition
\begin{align} 
{\cal C}^0(x)=c_0(x),
\nn 
\end{align}

From the fully discrete scheme \refe{c1}-\refe{c3} and the above weak formulation, 
we obtain the error equations
\begin{align} 
&
\left( D_{\tau}\xi_c^{n+1},\phi_h\right)
+  \left(D(\u_h^n)\nabla \xi_c^{n+1},\nabla \phi_h \right) + \left( \xi_c^{n+1} q^P, \phi_h \right)\ 
\label{s1} \\
&
=  \frac{1}{\tau} 
\left( e_c^{n+1} - e_c^n(x_{U^n}), \phi_h \right )
- \frac{1}{\tau}  \left(\wt \C^n_h(x_{U^n}) - \wt \C^n_h(x_{\u_h^n}),\phi_h\right) 
+ \frac{1}{\tau}  \left(\xi_c^n(x_{\u_h^n}) - \xi_c^n,\phi_h\right) 
\nn \\  
& \quad 
+\left((D(\wt U^{n}_h)-D(\u_h^{n}))\nabla \wt \C_h^{n+1}, \nabla\phi_h \right)
+ \left( e_c^{n+1}q^P, \phi_h \right) 
\qquad \phi_h\in V_h^1
\nn\\
&:=\sum_{i=1}^5 J_i,   
\nn  \\ 
&
\left( \frac{\mu(c_h^n)}{k(x)} \u_h^n
- \frac{\mu(\C^n)}{k(x)} U^n, \v_h \right)
=\left(p_h^n - P^n, \nabla\cdot \v_h \right), \qquad 
\v_h\in H_h^0,
\label{s2}\\
&
\left(\nabla\cdot (\u_h^n -U^n ),\varphi_h \right)=0, \qquad 
 \varphi_h\in S_h^0 \, . 
\label{s3} 
\end{align} 
where $\xi_c^0 = I_h c_0 - \widetilde c_0$. 

By \refe{e-proj-2} and \refe{s3}, we can see that
$(\nabla\cdot\xi_u^n,\varphi_h)=0$ for any $\varphi_h\in H_h^0$, 
which implies 
\begin{align} 
\nabla \cdot \xi_u^n = 0 \, . 
\label{div-u}
\end{align} 
and therefore, by noting the definition of $H^0_h$, at each element, $\xi_u^n$ is a constant vector.

To prove Theorem \ref{main}, first we rewrite \refe{s2} into
\begin{align}
\left( \frac{\mu(c^n_h)}{k(x)} \xi_u^n, \v_h \right ) 
+  \left ( \left ( \frac{\mu(c^n_h)}{k(x)} - \frac{\mu(\C^n)}{k(x)} \right ) 
\widetilde U^n_h, \v_h \right )    
=  
\left ( \xi_p^n, \nabla \cdot \v_h \right ) 
\label{p-2}
\end{align} 
where we have noted \refe{e-proj-1}.

By taking $\v_h = {\color{red}| \xi_u^n|^2  \xi_u^n} \in S_h^0$ in the last equation, we see that  
\begin{align} 
 \| \xi_u^n \|_{{\color{red} L^4}} 
&  \le \| \frac{\mu(c^n_h)}{k(x)} - \frac{\mu(\C^n)}{k(x)} \|_{L^4}
\| \widetilde U^n_h \|_{L^\infty}     
\nn \\ 
& \le 
 C \| c^n_h - \C^n \|_{L^4} 
\nn \\ 
& \le  C \| \xi_c^n \|_{{\color{red} H^1}}    + C h^2 
\label{xiu} 
\end{align} 
where we have noted \refe{mu} and \refe{U-infty}  and used Lemma \ref{l3-4}.

Secondly we prove the following 
primary estimate by mathematical induction 
\begin{align} 
\| \xi_{c}^n \|_{L^2} + \tau^{1/2} \| \nabla \xi_{\color{red} c}^n \|_{L^2}  
 \leq h^{11/6},    \quad n=0,1,\cdots,N   \, . 
\label{pri}  
\end{align}

Since $\xi_c^0=c_h^0-I_h\C^0=0$,
the estimate \refe{pri} holds for $n=0$.

We assume that \refe{pri} holds for
$n\leq m$ for some integer $m\geq 0$, which with \refe{U-infty}, 
\refe{chh}, \refe{xiu}, 
Lemma \ref{l3-3} and inverse 
inequalities implies 
\begin{align} 
& \| c_h^n \|_{L^\infty}  \le \| \wt {\C}_h^n \|_{L^\infty} + 
 \| \xi_c^n \|_{L^\infty}  
\le C + C h^{-3/2} \| \xi_c^n \|_{L^2} 
\le C 
\label{ch} 
\\ 
&  \| \u_h^n \|_{{\color{red} L^\infty}} \le \| \wt U_h^n \|_{L^\infty} + \| \xi_u^n \|_{L^\infty} 
\le C + C h^{-3/2} (\| \xi_c^n \|_{L^2} +h^2) 
\le C \, . 
\label{uh} 
\end{align} 
By  Lemma \ref{l3-3},  $U^n \in W^{2,3}(\Omega)\cap C^1(\Omega)$ satisfies 
the condition  \refe{Lip}. Taking $\phi_h=\xi_c^{n+1}$ in \refe{s1}, 
by Lemma \ref{l3-2} (ii)  with $p=2$, we have
\begin{align*} 
J_1(\xi_c^{n+1}) &=
\left(D_{\tau} e_c^{n+1}, \xi_c^{n+1}\right)
+ \frac{1}{\tau}  \left( e_c^n - e_c^n(x_{U^n}) , \xi_c^{n+1}\right) 
\nn\\
&\le
C\|D_{\tau} e_c^{n+1} \|_{L^2}  
\|\xi_c^{n+1}\|_{L^2}  
+C 
\| e_c^n \|_{L^2}  \|U^n\|_{W^{1,\infty}}\|\xi_c^{n+1}\|_{H^1} 
\nn \\  
& \le \epsilon \|\xi_c^{n+1}\|_{H^1}^2 
+ C\|D_{\tau} e_c^{n+1} \|_{L^2}^2  
+ C h^4   
\nn 
\end{align*} 
and
\begin{align*}  
J_5(\xi_c^{n+1}) &\leq
C\|q^p\|_{L^3}\| e_c^{n+1} \|_{L^2} 
\|\xi_c^{n+1}\|_{L^6} 
\nn\\
&\leq \epsilon\|\xi_c^{n+1}\|_{H^1}^2 + C_{\epsilon} h^4 \, .  
\end{align*} 
By noting \refe{U-infty}, \refe{uh} and \refe{chh}, we have 
\begin{align} 
\| (D(\wt U_h^{n}) - D(\u_h^{n}) ) \nabla \widetilde \C_h^{n+1} \|_{L^2} 
& = \|D'(\chi) \xi_u^{n} \cdot  \nabla \widetilde \C_h^{n+1}\|_{L^2}  
  \le C \| \xi_u^{n} \|_{L^2} 
  \nn 
  \end{align} 
and therefore, 
\begin{align*} 
J_4(\xi_c^{n+1}) \leq
C\|\xi_u^n\|_{L^2}  
\|\nabla\xi_c^{n+1}\|_{L^2}  \, . 
\nn
\end{align*} 
Moreover,  by \refe{U-infty} and \refe{pri},  
\begin{align}  
\| U^n - \u_h^n \|_{L^\infty} & 
\le \| \xi_u^n \|_{L^\infty} + \| e_u \|_{L^\infty} 
\nn \\ 
& \le C h^{-3/4} \| \xi^n_u \|_{L^4} +  C h \log(1/h)  
\nn \\ 
& \le C h^{-3/4} \| \xi^n_{\color{red} c } \|_{H^1} +  C h \log(1/h)  
\nn \\ 
& \le Ch^{13/12} \tau^{-1/2} +  C h \log(1/h)  
\nn 
\end{align} 
and 
\begin{align} 
\tau \| U^n - \u_h^n \|_{L^\infty} 
\le C h^{13/12}\tau^{1/2} + Ch \tau \log(1/h)\, . 
\end{align} 
Then both $U^n$ and $U^n-\u_h^n$ satisfy  the condition \refe{Lip}. 
By Lemma \ref{l3-2} (i) and (iii), we can see that 
\begin{align} 
J_2(\xi_c^{n+1}) &= 
\frac{1}{\tau}  \left( \C^n(x_{U^n}) - \C^n(x_{\u_h^n}), \xi_c^{n+1} \right) 
-  \frac{1}{\tau}  \left(  e_c^n(x_{U^n}) -  e_c^n(x_{\u_h^n}), \xi_u^{n+1}  \right) 
\nn \\ 
 &\le   C \| \C^n \|_{W^{2,4}}  (\| U^n -  \u_h^n \|_{H^{-1}}  + \tau \| U^n -  \u_h^n \|_{L^4}^2 ) 
  \| \xi_c^{n+1} \|_{H^1} 
 \nn \\ 
 &  + C\| e_c^n \|_{W^{1,3}} \| U^n - \u_h^n \|_{L^2} \| \xi_c^{n+1} \|_{L^6}  
 \nn 
 \end{align} 
 By \refe{xiu}, inverse inequalities and mathematics induction, 
\begin{align} 
 \| U^n -  \u_h^n \|_{L^4} & \le C (\| e_u^n \|_{L^4} + 
 \| \xi_u^n \|_{L^4} ) 
 \nn \\ 
 &  \le C h + 
C h^{-3/4} \| \xi_u^n \|_{L^2} 
\nn \\ 
& \le C h + 
C h^{-3/4} ( \| \xi_c^n \|_{L^2} + h^2) 
\nn \\ 
& \le  C h 
\nn 
\end{align} 
and therefore, by \refe{Up-H-1} 
\begin{align} 
J_2(\xi_c^{n+1}) & \le   
C \left ( \| \xi_u^n \|_{L^2} + \| e_u^n \|_{H^{-1}} + \tau h^2 
+ \| e_c^n \|_{W^{1,3}} (\| \xi_u^n \|_{L^2} + \| e_u^n \|_{L^2}) \right ) 
 \| \xi_c^{n+1} \|_{H^1} 
\nn \\ 
 & \le 
C ( \| \xi_c^n \|_{L^2} +  h^2  )  \| \xi_c^{n+1} \|_{H^1}  
\nn \\ 
 & \le 
\epsilon \| \xi_c^{n+1} \|_{H^1}^2 + 
C  (\| \xi_c^n \|_{L^2}^2 + h^4 )  \, . 
 \nn 
\end{align} 

For $J_3$, by Lemma \ref{l3-2} (i)-(ii),  we have 
\begin{align} 
|J_3(\xi_c^{n+1}) |&= 
\frac{1}{\tau} |(\xi_c^n(x_{\u_h^n})-\xi_c^n(x_{U^n}), \xi_c^{n+1})| 
+ | \frac{1}{\tau} (\xi_c^n(x_{U^n}) - \xi_c^n \, , \xi_c^{n+1}) | 
\\
&\leq
C\|\xi_c^n\|_{H^1}
\|\u^n_h - U^n \|_{L^3} \|\xi_c^{n+1}\|_{L^6}
+ C\|\xi_c^n\|_{L^2} \|U^n\|_{W^{1,\infty}}\|\xi_c^{n+1}\|_{H^1} 
\nn \\ 
& \le
\epsilon (\| \xi_c^{n+1} \|_{H^1}^2 + \|\xi_c^n\|_{H^1}^2)  
+ C  \| \xi_c^{n} \|_{L^2}^2 
\nn  
\end{align} 

Substituting the above estimates into 
\refe{s1}, we obtain 
\begin{align} 
& \frac{\|\xi_c^{n+1}\|_{L^2}^2-\|\xi_c^n\|_{L^2}^2}{\tau}
+  \left\|\sqrt{D(\u_h^n)}\nabla\xi_c^{n+1}\right\|_{L^2}^2 
\nn\\
& \quad \leq
C(\|\xi_c^{n+1}\|_{L^2}^2+\|\xi_c^n\|_{L^2}^2)
+\epsilon(\|\nabla\xi_c^{n+1}\|_{L^2}^2
+\|\nabla\xi_c^n \|_{L^2}^2)+C h^4 
\nn\\
& \quad +C\|D_{\tau} e_c^{n+1} \|_{L^2}^2
\label{en-1} 
\end{align} 

By Gronwall's inequality, we arrive at 
\begin{align} 
& \|\xi_c^{n+1}\|_{L^2}^2
+\sum_{m=0}^n\tau\|\nabla\xi_c^{m+1}\|_{L^2}^2 \leq Ch^4 + C \sum_{m=0}^n \tau \| 
D_{\tau} e_c^{m+1} |_{L^2}^2
 \le h^{11/3} 
\nn 
\end{align} 
when $h \le h_1$ and  $\tau \le \tau_1$ for some $\tau_1>0$, 
where we have noted 
\refe{p2-2}. 
Therefore, the induction is closed and \refe{pri} holds for any $n \ge 0$. Moreover, we have 
\begin{align} 
\| \xi_c^{n+1} \|_{L^2} \le C h^2  
\label{err-1} 
\end{align} 
and by \refe{div-u} and \refe{xiu}, 
\begin{align} 
\| \xi_u^{n+1} \|_{H(div)} \le C h^2 
\label{xiuh} 
\end{align} 
which with the error estimate \refe{Up-Lp-div}  for mixed projection and the error estimate \refe{p2-1} 
for the quasi-projection leads to 
\begin{align} 
\| c_h^n - {\cal C}^n \|_{L^2} + h \| \u_h^n - U^n \|_{H(div)} \le C h^2, \qquad n=1,2,...,N \, . 
\label{cu}
\end{align}

To derive an estimate for $\| p_h^n - P^n \|_{L^2}$,  we follow a traditional way 
used in \cite{DEW-2, DEW, ERW}. 
From \refe{e-proj-1}-\refe{e-proj-2} and \refe{s2}-\refe{s3}, we see that 
\begin{align} 
& \left ( \frac{\mu(c_h^n)}{k(x)} (\u_h^n - \widetilde U^n_h), \v_h \right ) =  (p_h^n - \widetilde P^n_h, \,  
\nabla \cdot \v_h) - \left ( \frac{\mu(c_h^n) - \mu(C^n)}{k(x)} \widetilde U^n_h, \v_h \right ), \quad 
\v_h \in H_h^0
\\ 
& \, ( \nabla \cdot ( \u_h^n - \widetilde U^n_h ), \varphi_h ) = 0, \qquad 
 \varphi_h\in S_h^0 \, . 
\end{align} 
By Brezzi's Proposition 2.1 in \cite{Brezzi},  the error of the pressure  is bounded by 
\begin{align} 
\| p_h^n - \widetilde P^n_h \|_{L^2} \le C (1+ \| \widetilde U_h^n \|_{L^\infty} ) \| c_h^n - C^n \|_{L^2} 
\, . 
\end{align} 
By noting the bound \refe{uh} for $\widetilde U^n$, the projection error estimate \refe{Up-Lp-div} 
and \refe{cu}, we get 
\begin{align*}
\|p_h^n-P^n\|_{L^2}\leq Ch
\quad n=0,1,\cdots, N \, . 
\end{align*}
Finally, 
\refe{error0}-\refe{error} follow \refe{cu}, the last equation and 
Lemma \ref{l3-3}. The proof is complete. \quad 
\endproof
\mbox{}\vskip0.1in 

For the upgraded numerical pressure/velocity $(\widehat p_h^m, \widehat \u_h^m)$ generated 
by the post-process \refe{new-p}-\refe{new-u}, 
by taking  a similar approach, we can see that 
\begin{align} 
\| \xi_{\widehat u} \|_{L^2} \le \| \xi_c \|_{L^2} + C h^2 \le C h^2 
\nn 
\end{align} 
and 
\begin{align} 
\| \widehat p_h^m - P^m \|_{L^2} \le C  h^2 \, . 
\nn 
\end{align} 
The corollary \ref{u-Error} follows immediately. 
\quad \endproof 
\mbox{}\vskip0.1in 

{\it Remarks.} In some practical cases, one may use two different partitions, 
$\pi_{h_p}$ and $\pi_{h_c}$ 
for the concentration equation and pressure/velocity equation, respectively. 
Previous analysis based on  two different meshes requires  
certain mesh condition, which excluded the most commonly-used mesh $\pi_{h_p}= \pi_{h_c}$. 
It is possible to extend our approach to the problem with two different partitions to establish 
 the general error estimate 
\begin{align} 
& \| c_h^m - c^m \|_{L^2} \le C_0 ( \tau + h_p^2 + h_c^2 ) 
\\ 
& \| \u_h^m - \u^m \|_{L^2} + \| p_h^m - p^m \|_{L^2} 
\le C_0 ( \tau + h_p + h_c ) 
\end{align} 
 under the condition $h_c \ge Ch_p$  for any given $C>0$, which includes the case $h=h_p=h_c$. 
Also the extension to the general finite element spaces $(V_h^r, S_h^k, H_h^k)$ is possible, 
while 
higher regularity of the solution of the system is required. 

\section{Numerical results}
\setcounter{equation}{0} 
In this section, we present some numerical results 
in both two and three-dimensional 
porous media to confirm our theoretical analysis and show the efficiency of the post-processing. 
We always assume that the solution of the system is smooth. 
The problem with non-smooth solutions was considered in \cite{CLLS,LS3}. 
 Computations are performed by the free software FreeFem++ 
 for two-dimensional case and FEniCS for three-dimensional case.

We rewrite the system (\ref{e1})-(\ref{e3}) by
\begin{align} 
  & \frac{\partial c}{\partial t} - \nabla \cdot \left( D({\u})\nabla c \right)
  +{\u}\cdot\nabla c = g,
  \label{n1}  \\
  & \u = -\frac{1}{\mu(c)} \nabla p,
  \label{n2}  \\
  & \nabla \cdot \u = f \, .
  \label{n3}
\end{align}

\begin{figure}[ht]
  \centering
  \begin{tabular}{cc}
    \includegraphics[width = 50mm]{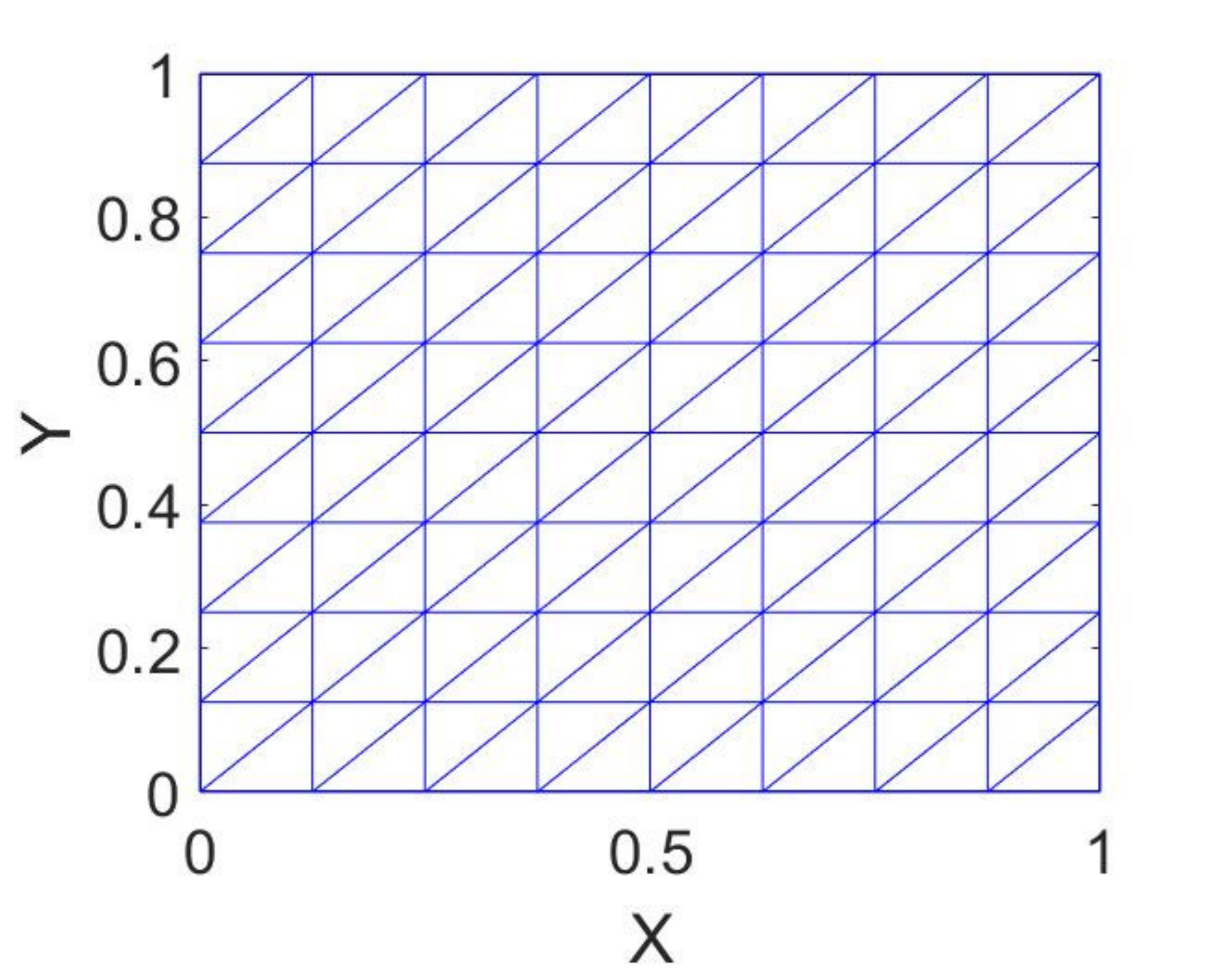}
    & 
    \epsfig{file=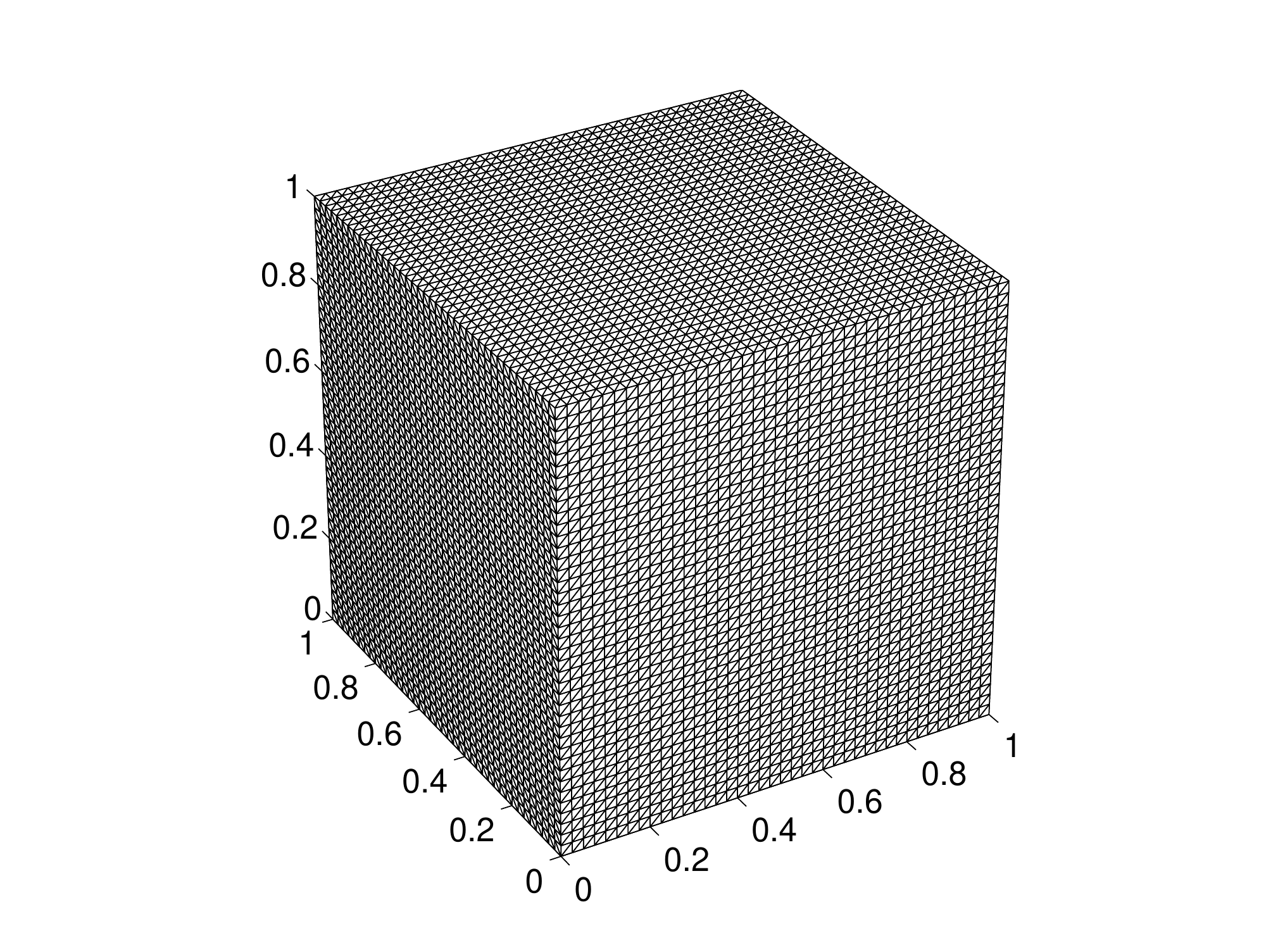,width=60mm}
  \end{tabular}
  \caption{A uniform triangular/tetrahedron mesh on the unit square/cube.}
  \label{mesh-2d}
\end{figure}

First, we consider the system \refe{n1}-\refe{n3} 
on the unit square domain $\Omega = [0,1] \times [0,1]$, where 
$D(\mathbf{u})=\frac{1}{40}(1+|\mathbf{u}|^2)$
and $\mu(c)=1/(1+c^2)$.
We set the terminal time $T = 1.0$. The functions $f$ and $g$
are chosen correspondingly to the exact solution
\begin{eqnarray}
&&c=1+20e^t(1+t^2)\sin(x^2)\sin(y^2)(1-x)^2(1-y)^2,\\
&&p=3+400e^t(1+t^3)x^2y^2(1-x)^3(1-y)^3. 
\end{eqnarray}
Numerical simulations of mixed FEM and characteristics mixed FEM with 
the FE space $(V^1_h, S_h^1, H_h^1)$ were presented in \cite{LS1} and \cite{WSS}, 
respectively. 

Here, we solve the system
\refe{n1}-\refe{n3} by the numerical scheme \refe{c1}-\refe{c3} with the lowest-order 
characteristics mixed FEM,
in which a uniform triangular partition with $M+1$ nodes in each
direction is used,  see Figure \ref{mesh-2d} for an illustration with $M=8$, where $h=\sqrt{2}/M$.
To show the optimal convergence rates, we choose $\tau=h^2$ in our numerical simulations. 
We present in Table 1 numerical results at $t=1.0$. From 
Table 1, 
we can observe clearly the second-order convergence rate for the concentration in $L^2$-norm, 
which is the most important physical component in applications. 
The convergence rate for both the pressure and Darcy velocity is $O(h)$ in $L^2$-norm, which is optimal in the traditional sense.  The one-order lower approximation to 
$(p, u)$ does not affect the accuracy of the numerical concentration.  

On the other hand, we resolve the system \refe{new-p}-\refe{new-u} with the obtained 
concentration $c_h^N$ for $(\widehat p_h^N,  \widehat \u_h^N ) \in (S_h^1, H_h^1)$. 
We present these numerical results in Table 2 from which 
we can see  the second-order accuracy for both pressure and Darcy velocity. For comparison, 
we present in Table 3 numerical results of the scheme 
\refe{c1}-\refe{c3} with $(c_h^n, p_h^n, \u_h^n ) \in (V^1_h, S_h^1, H_h^1)$. 
Numerical results show that the post-processing defined in \refe{new-p}-\refe{new-u} provides 
the same accuracy as the classical FEM solution 
in  $ (c_h^n, p_h^n, \u_h^n) \in (V_h^1, S_h^1, H_h^1), n=1,2,...,N$, while the former requires much less 
computational cost. 

To test the stability of the scheme, we solve the system \refe{new-p}-\refe{new-u} with several different $M$ for each $\tau=1/20, 1/30, 1/40$. Numerical results presented in Figure 2 illustrates 
 that 
 $L^2$-norm errors converge to $O(\tau)$ as $M$ increases. 
This shows that the time step restrictions given in those previous works are not necessary.

\begin{table}[h]
\caption{$L^2$ errors of the scheme \refe{c1}-\refe{c3} with $(c_h^n, p_h^n, \u_h^n ) \in (V^1_h, S_h^0, H_h^0)$ 
in 2D. 
}
\label{order}
\begin{center}
\begin{tabular}{c|c|c|cc}
  \hline
 $\tau=1/M^2$ & $\|c(\cdot,t_N)-c_h^N\|_{L^2}$  &
 $\|\mathbf{u}(\cdot,t_N)-\u_h^N\|_{L^2}$  &  $\|p(\cdot,t_N)-p_h^N\|_{L^2}$   \\ \hline 
  $M=8$   & 3.923e-2 & 5.945e-1 & 1.532e-1  \\ 
  $M=16$	& 7.518e-3 &3.020e-1 & 7.763e-2  \\ 
  $M=32$	& 1.704e-3 & 1.517e-1& 3.759e-2   \\ 
  $M=64$  & 4.250e-4 & 7.681e-2 & 1.882e-1 \\ \hline 
  Convergence  order &   2.00  & 0.98       & 1.00                 \\
   \hline
\end{tabular}
\end{center}
\end{table}

\begin{table}[h]
\caption{$L^2$-norm errors of the post-processing  with $( p_h^N, \u_h^N ) \in (S_h^1, H_h^1)$ in 2D. 
}
\label{order}
\begin{center}
\begin{tabular}{c|c|cc}
  \hline
 $\tau=1/M^2$ &  
 $\|\mathbf{u}(\cdot,t_N)-\widehat \u_h^N\|_{L^2}$  &  $\|p(\cdot,t_N)-\widehat p_h^N\|_{L^2}$  
  \\ \hline 
  $M=8$   & 3.299e-2 & 8.639e-2 \\
  $M=16$	& 7.025e-2 & 2.255e-2  \\
  $M=32$	& 1.588e-3 & 5.607e-3    \\ 
  $M=64$ &  3.871e-4 & 1.401e-3 \\  \hline 
   order &   2.04     & 2.00     \\ 
   \hline
\end{tabular}
\end{center}
\end{table}

\begin{table}[h]
\caption{$L^2$-norm errors of scheme  \refe{c1}-\refe{c3} with $(c_h^n, p_h^n, \u_h^n ) \in (V_h^1, S_h^1, H_h^1)$ 
in 2D.}
\label{order}
\begin{center}
\begin{tabular}{c|c|c}
  \hline
 $\tau=1/M^2$ &     
 $\|\mathbf{u}(\cdot,t_N)- \widehat \u_h^N\|_{L^2}$  &  $\|p(\cdot,t_N)- \widehat p_h^N\|_{L^2}$   
  \\ \hline 
  $M=8$   &  3.791e-2 & 8.625e-2\\
  $M=16$	&  7.901e-3 & 2.221e-2  \\
  $M=32$	& 1.801e-3 & 5.601e-3   \\ 
  $M=64$  & 4.451e-4 & 1.423e-3 \\ \hline 
   order &    2.01  & 1.98    \\ 
   \hline
\end{tabular}
\end{center}
\end{table}

\begin{figure}[h]
\begin{center}
\begin{tabular}{c}
\epsfig{file=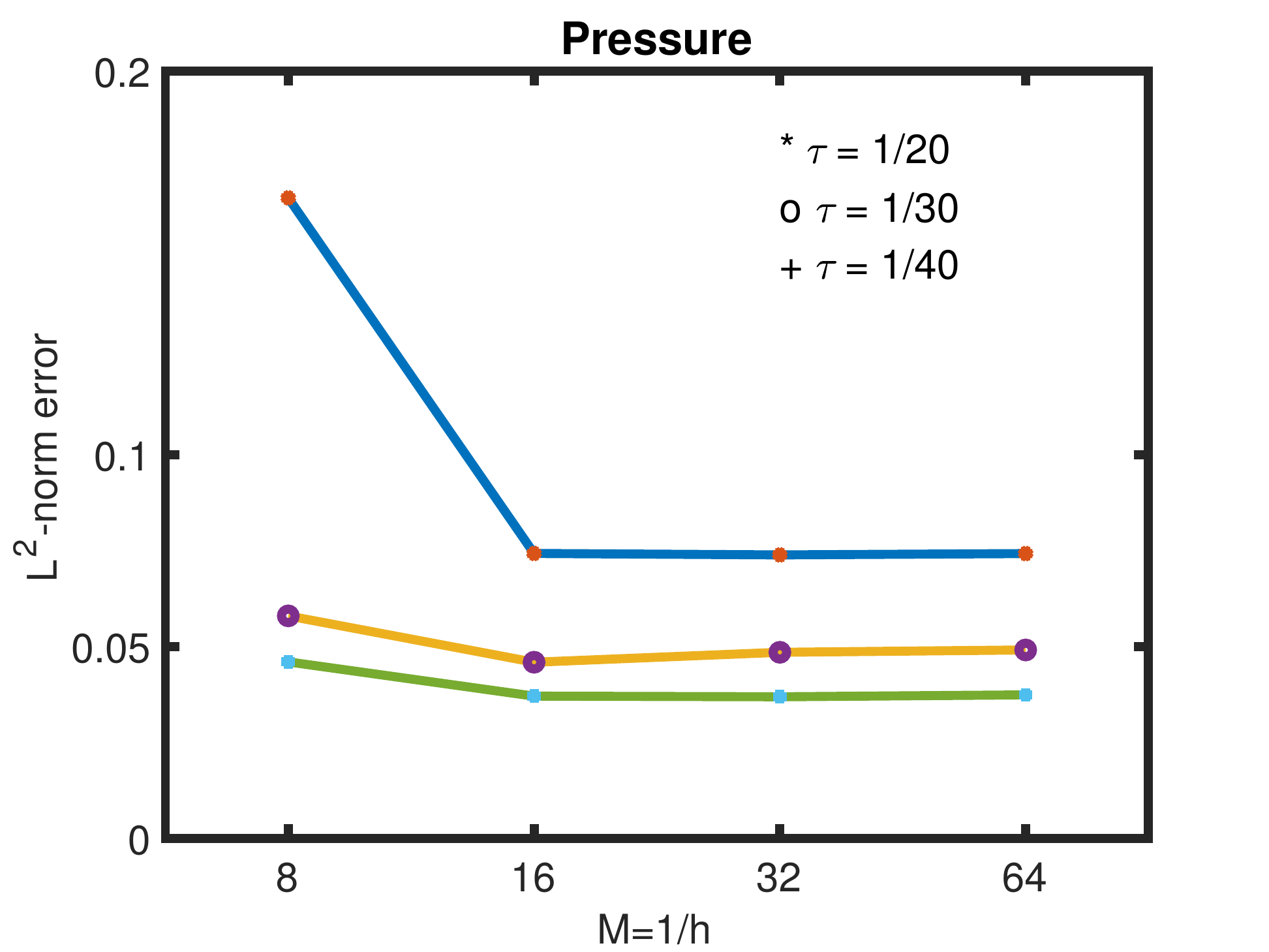,height=1.8in,width=2.7in} \quad
\epsfig{file=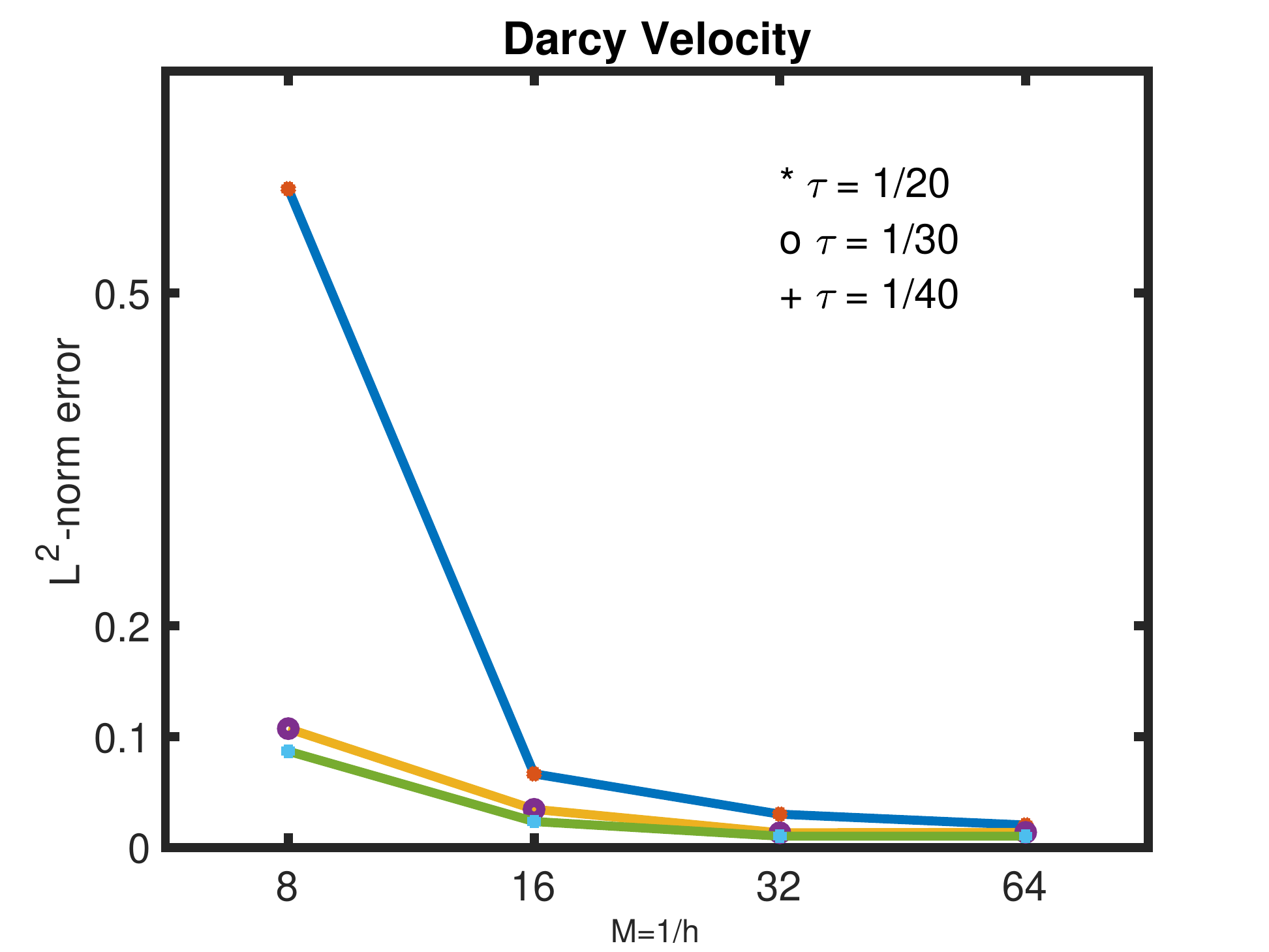,height=1.8in,width=2.7in} \\ 
\\ 
\epsfig{file=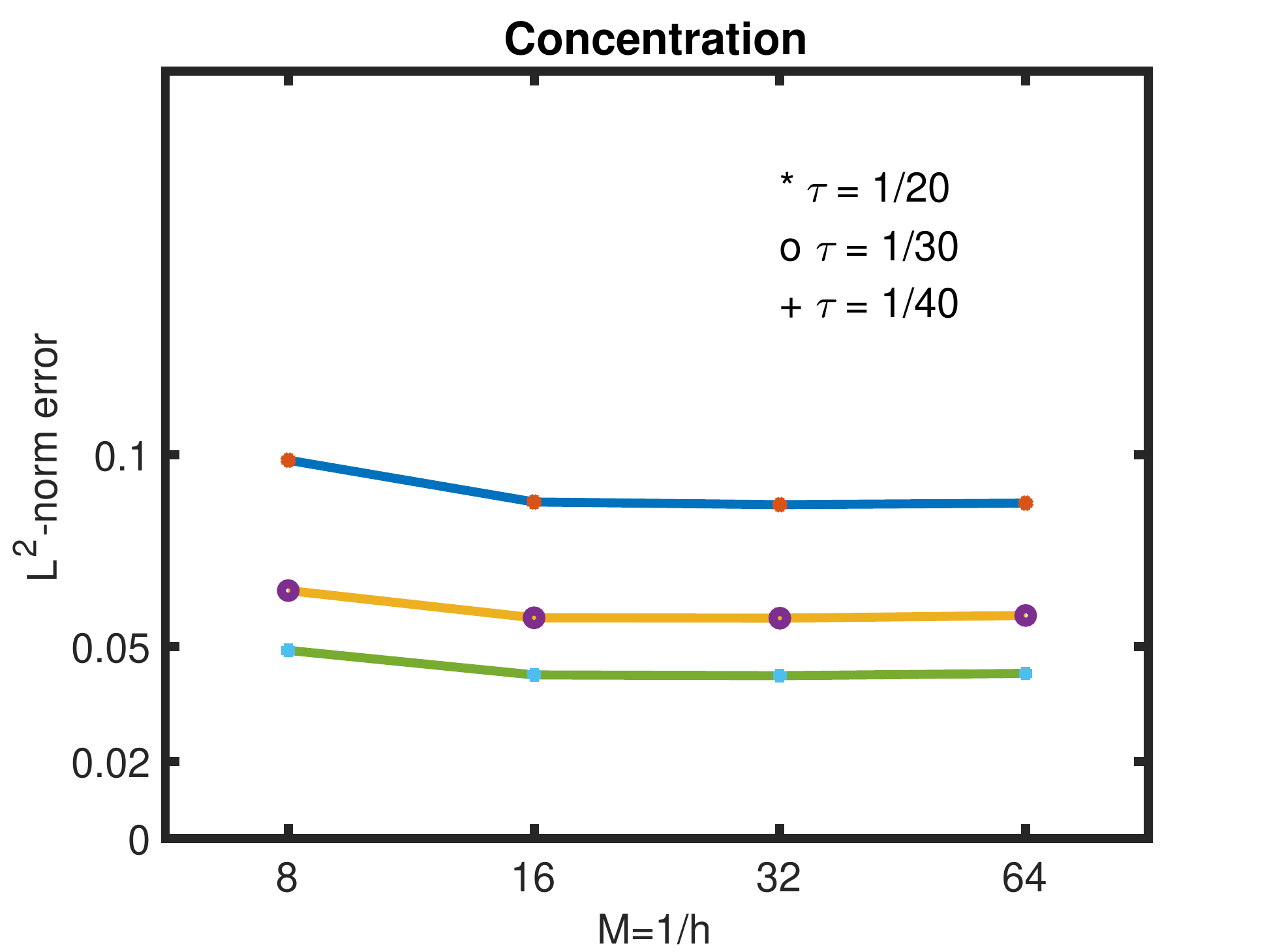,height=1.8in,width=2.7in}
\end{tabular}
\end{center}
\caption{Stability of the scheme \refe{c1}-\refe{c3}}
\end{figure}

Secondly, we study the system \refe{n1}-\refe{n3}
in a three-dimensional cube $[0,1] \times [0,1] \times [0,1]$, 
where $D(\u)=1+\u\otimes\u$ and 
$\mu(c)=1+c^2$ and 
the functions $f$ and $g$ are chosen correspondingly 
to the smooth exact solution
\begin{align}
  & c =  \frac{1}{10} \exp(-t) \sin^2(\pi x) \sin(\pi y) \sin(\pi z) \, ,
 \label{num003}\\
 & p = \exp(-t) \cos( \pi x) \cos( \pi y) \cos( \pi z) \, .
 \label{num004}
\end{align}
We also set the terminal time T = 1.0 in this example.

We use a uniform tetrahedra mesh with $M+1$ nodes 
in each direction
($h = \frac{\sqrt{2}}{M}$), see Figure \ref{mesh-2d}. We solve the system 
\refe{n1}-\refe{n3} on the unit cube with $\tau=1/512$ 
and $M=8,16,32$. We present our numerical results in 
Table \ref{linear-L2-1}. Numerical results confirm 
the second-order 
accuracy of concentration by the lowest characteristic-mixed 
FEM.

Again, after getting $c_h^N$, we resolve \refe{n2}-\refe{n3}
at the terminal time $T=1.0$ with $(\widehat{\u}_h^{N},\widehat{p}_h^{N}) \in (S_h^1, H_h^1)$. 
We present the $L^2$-norm errors of the recovered
numerical solution $(\widehat{\u}_h^{N},\widehat{p}_h^{N})$
in Table \ref{linear-L2-1-recover}.
The second-order accuracy of numerical solution 
$(\widehat{\u}_h^{N},\widehat{p}_h^{N})$
is observed clearly, which confirms that the approximation for $(\u,p)$ in three dimensions can also be significantly improved by the proposed 
post-processing.

\section{Conclusion}
\setcounter{equation}{0}
We have established  optimal error estimates of the commonly-used lowest-order characteristics-mixed FEMs with linearized Euler scheme for miscible displacement problems 
under a weak time step condition 
$\tau = o(1/|\log h|)$. Previous analysis only provided a sub-optimal estimate for the concentration. 
We have shown theoretically and numerically  that the lower-order approximation to 
the velocity/pressure does not pollute the numerical concentration 
and also, the scheme allows one to use a large time step. 
The analysis presented in this paper can be easily extended to 
other existing methods, such as ELLAM and high-order characteristic approximations. The analysis presented in this paper is based on the assumption of certain strong regularity of the solution. 
The problem with weaker regularity assumption is of interest. 
Some existing works can be found in literature, such as 
\cite{AH, DPBL} for mixed finite volume methods and \cite{CDL,DEPT}  
for a framework of gradient discretization methods, including 
mixed FE-ELLAM and hybrid mimetic mixed-ELLAM schemes.  
On the other hand, 
 theoretical analysis in this paper
is based on the $\Omega$-periodic model as usual \cite{DRW,Dur,ERW,Rus,WSS}
to avoid the technical difficulties on the boundary. This periodic assumption is
physically reasonable. For the problem with Neumann boundary conditions, some further approximation to $c_h^n(x)$ 
was mentioned in \cite{Rus}.

\begin{table}[h]
  \caption{$L^2$-norm errors of the scheme \refe{c1}-\refe{c3} 
   with $(c_h^n, p_h^n, \u_h^n ) \in (V_h^1, S_h^0, H_h^0)$ 
  in 3D.}
     \begin{center} 
   \label{linear-L2-1}
    \begin{tabular}{c|cccc}
      \hline
      \hline
      \rule{0pt}{3ex}
      $\tau = 1/512$ & $c_h^{N}-c^N\|_{L^2}$ &
      $\|\u_h^{N}-\u^N\|_{L^2}$ & $\|p_h^{N}-p^N\|_{L^2}$
      \rule[-1.2ex]{0pt}{0pt} \\\hline      
      M = 8  	& 2.441e-03  & 2.512e-01 & 4.872e-02 \\ 
      M = 16  	& 6.544e-03 & 1.263e-01 & 2.451e-02 \\ 
      M = 32  	& 4.422e-04  & 6.284e-02 & 1.221e-02 \\ \hline      
      Order 	& 1.94        & 1.00        & 1.00       \\\hline
      \hline 
    \end{tabular}
  \end{center}
\end{table}

\begin{table}[h]
   \caption{$L^2$-norm errors of the post-processing  
with $( \widehat p_h^N, \widehat \u_h^N ) \in (S_h^1, H_h^1)$ 
 in 3D .}
    \begin{center} 
   \label{linear-L2-1-recover}
   \begin{tabular}{c|cc}
     \hline
    \hline
     \rule{0pt}{3ex}
    $\tau = \frac{1}{2 M}$ & $\|\widehat \u_h^{N}-\u^N\|_{L^2}$
    & $\|\widehat p_h^{N}-p^N\|_{L^2}$
    \rule[-1.2ex]{0pt}{0pt} \\\hline
     M = 8 	& 1.591e-02  & 1.221e-02 \\ 
     M = 16  	& 4.143e-03  & 1.012e-02 \\ 
      M = 32 	& 1.046e-03  & 1.534e-03 \\ \hline    
      Order 	& 1.99        & 1.99       \\\hline
      \hline 
    \end{tabular}
  \end{center}
\end{table}
\vskip0.1in

\noindent{\bf Acknowledgments}~The author would like to thank the anonymous referee
for the careful review and valuable suggestions and comments, which have greatly  
improved this article.


\begin{thebibliography}{99}


\bibitem{AH}
T.  Arbogast and C.S.~Huang, A fully mass and volume conserving implementation of a characteristic method for transport problems, 
{\em  SIAM J. Sci. Comput.}, 28(2006), 2001--2022. 
 
\bibitem{ASW} M. Al-Lawatia, R.C. Sharpley and H. Wang,
Second-order characteristic methods for advection-diffusion equations
and comparison to other schemes, {\em Adv. Water. Resour.}, 22 (1999), 741-768.

\bibitem{BB}
J.~Bear and Y.~Bachmat,
{\em Introduction to Modeling of Transport Phenomena in Porous Media},
Springer-Verlag, New York, 1990.

 \bibitem{BNV} 
 A.~Bermudez, M.R.~Nogueiras and C.~Vazquez, 
 Numerical analysis of convection-diffusion-reaction problems with higher order characteristics/finite elements. II. Fully discretized scheme and quadrature formulas, 
 {\em  SIAM J. Numer. Anal.}, 44 (2006), 1854-1876.
 
 \bibitem{BS}
    S.~Brenner and L.~Scott,
    {\em The Mathematical Theory of Finite Element Methods},
    Springer, New York, 2002.


 \bibitem{Brezzi} 
 F.~Brezzi, On the existence, uniquness and appproximation of saddle-point problems, arising from Lagrangian multipliers, {\em RAIRO
Anal. Numer.}, 2 (1974), 129--151. 


 \bibitem{CLLS}
W.~Cai, B~ Li, Y.~Lin  and W.~Sun, 
Analysis of fully discrete FEMs for miscible displacement
in porous media with Bear-Scheidegger diffusion-disperson tensor.
{\em Numer Math},  141 (2019), 1009--1042. 

\bibitem{CRHE} 
M.A.~Celia, T.F. Russell, I. Herrera and R.E. Ewing,
An Eulerian-Lagrangian localized adjoint method for the advection-diffusion equation,
{\em Adv. Water. Resour.}, 13 (1990), 187--206.

\bibitem{CCW}
F. Chen, H. Chen and H. Wang,  
An optimal-order error estimate for a Galerkin-mixed finite element time-stepping procedure for porous media flows, 
{\em Numer. Methods Partial Differ. Equations}, 28.2(2012), 707--719.

\bibitem{CWW}
 A. Cheng, K. Wang and H. Wang,  
Superconvergence for a time-discretization procedure for the mixed finite element 
approximation of miscible displacement in porous media,  
{\em Numer. Methods Partial Differ. Equations}, 28(2012), 1382--1398.


\bibitem{CDL}
 H. M.~Cheng,  J.~Droniou and K.N.~Le, 
 Convergence analysis of a family of ELLAM schemes for a fully coupled model of miscible displacement in porous media, 
 {\em  Numer. Math.},  141(2019), 353--397.

 
\bibitem{DRW}
 C.N. Dawson, T.F. Russell and M.F. Wheeler,
Some improved error estimates for the modified method of
characteristics,
{\em SIAM J. Numer. Anal.}, 26 (1989), 1487-1512.


\bibitem{DPBL}
M.~ D'Elia, M. Perego, P.~Bochev and D.~Littlewood, 
A coupling strategy for nonlocal and local diffusion models with mixed volume constraints and boundary conditions, 
{\em  Comput. Math. Appl.}, 71 (2016), 2218--2230.

 
\bibitem{DO}
L.~Demkowicz and J.T.~Oden,
An adaptive characteristic Petrov-Galerkin finite element method
for convection-dominated linear and nonlinear parabolic problems in one space variable,
{\em J. Comput. Phys.}, 67 (1986), 188-213.


\bibitem{DEW-2}
J.~Douglas,~Jr., R.E.~Ewing and M.F.~Wheeler,
The approximation of the pressure by a mixed method in the simulation of miscible displacement 
{\em RAIRO
Anal. Numer.}, 17 (1983), 17--33.


\bibitem{DEW}
J.~Douglas,~Jr., R.E.~Ewing and M.F.~Wheeler,
A time-discretization procedure for a
mixed finite element approximation of miscible displacement in porous media, {\em RAIRO
Anal. Numer.}, 17 (1983), 249-265.

\bibitem{DR}
J.~Douglas, Jr. and T.F.~Russell,
Numerical methods for convection-dominated diffusion problems based on
combining the method of characteristics with finite element or
finite difference procedures,
{\em SIAM J. Numer. Anal.}, 19 (1982), 871-885.


\bibitem{DRob}
J. Douglas, JR., and J. E. Roberts, 
Global estimates for mixed methods for second order elliptic equations, 
{\em Math. Comput.}, 44(1985), 39--52.


\bibitem{DEPT}
J.~ Droniou, R.~Eymard, A.~Prignet and K. S.~Talbot,  
Unified convergence analysis of numerical schemes for a miscible displacement problem, 
{\em  Found. Comput. Math.},  19(2019), 333--374. 

 
\bibitem{Dur}
R.G.~Duran,
On the approximation of miscible displacement in porous media by a method
of characteristics combined with a mixed method,
{\em SIAM J. Numer. Anal.}, 25 (1988), 989-1001.


\bibitem{Dur2}
R.G.~Duran,
Error analysis in $L^p$, $1 \le p \le \infty$, for mixed finite element methods for linear and 
quasi-linear elliptic problems, 
{\em RAIRO Mod. Math. Anal. Numer.} , 22(1988), 371--387.


\bibitem{Ewing}
R. E. Ewing, ed, The mathematics of Reservoir Simulation, Frontiers in Applied Mathematics, SIAM, Philadelphia, PA, 1983.

\bibitem{ERW}
R.E.~Ewing, T.F.~Russell and M.F.~Wheeler, Convergence analysis of an approximation
of miscible displacement in porous media by mixed finite elements and a modified method
of characteristics, {\em Comput. Methods Appl. Mech. Engrg.}, 47 (1984), 73-92.

\bibitem{EWang}
R. E. Ewing and H. Wang,
A summary of numerical methods for time-dependent advection-dominated
partial differential equations,
{\em J. Comput. Appl. Math.}, 128 (2001), 423-445.

\bibitem{EW}
R.E.~Ewing and M.F.~Wheeler, Galerkin methods for miscible displacement problems in
porous media, {\em SIAM J. Numer. Anal.}, 17 (1980), 351-365.

\bibitem{Feng}
X. Feng, On existence and uniqueness results for a coupled system modeling miscible 
displacement in porous media, 
{\em J. Math. Anal. Appl.}, 194 (1995), 883--910.

\bibitem{FN}
X.~Feng and M.~Neilan,
A modified characteristic finite element method for a fully nonlinear
formulation of the semigeostrophic flow equations,
{\em SIAM J. Numer. Anal.},  47 (2009), 2952-2981.

\bibitem{GPP}
 A.O. Garder, D.W. Peaceman and A.L. Pozzi,
Numerical calculations of multidimensional
miscible displacement by the method of characteristics,
{\em Soc. Pet. Eng. J.}, 4 (1964), 26-36.

\bibitem{GS}
 H. Gao and W. Sun, 
Optimal error analysis of Crank-Nicolson lowest-order Galerkin-mixed FEM for incompressible miscible flow in porous media, 
{\em Numerical Methods for PDEs}, 36(2020), 1773--1789.  

\bibitem{GN}
    L. Gastaldi and  R. H. Nochetto, 
Sharp maximum norm error estimates for general mixed finite element approximations to second order elliptic equations, 
{\em M2AN}, 23 (1989), 103--128. 

\bibitem{KM}
J.~Kacur and M.S.~Mahmood, 
Solution of solute transport in unsaturated porous media by the method of characteristics, 
{\em  Numer. Methods Partial Differential Equations}, 19(2003), 732-761.
 
 
\bibitem{KHR}
S. V. Krishnamachari, L. J. Hayes and T. F. Russell,
A finite element alternating-direction method combined with a modified
method of characteristics for convection-diffusion problems,
{\em SIAM J. Numer. Anal.}, 26 (1989), 1462-1473.


 \bibitem{KY}
  S.~Kumar and S.~Yadav,
Modified method of
characteristics combined with finite volume element methods for  incompressible miscible displacement problems in porous media,
{\em Int. J. Partial. Differ. Equ.}, 2014.

\bibitem{LS1}
B.~Li and W.~Sun, 
Error analysis of linearized semi-implicit Galerkin finite element methods
for nonlinear parabolic equations, {\em Int. J. Numer. Anal. Model.}, 10 (2013), 622-633.

\bibitem{LS2}
B.~Li and W.~Sun,
Unconditional convergence and optimal error estimates of a
Galerkin-mixed FEM for incompressible miscible flow in porous media,
{\em SIAM J. Numer. Anal.}, 51 (2013), 1959-1977.

\bibitem{LS3}
B. Li and W. Sun, Regularity of the diffusion-dispersion tensor and error analysis of 
FEMs for a porous media flow, 
{\em SIAM J. Numer. Anal.}, 53(2015), 1418--1437.


\bibitem{LWS}
 B. Li, J. Wang and W. Sun,
The stability and convergence of fully discrete Galerkin FEMs for
incompressible miscible flows in porous media,
{\em Commun. Comput. Phys.}, 15(2014), 1141--1158.

\bibitem{LWC}
D.~Liang, W.~Wang and Y.~Cheng,
An efficient second-order characteristic finite element method
for non-linear aerosol dynamic equations,
{\em Int. J. Numer. Methods Engrg.}, 80 (2009), 338-354.

\bibitem{MLY}
N.~Ma, T.~Lu and D.~Yang, Analysis of incompressible miscible displacement in porous media
by characteristics collocation method, {\em Numer. Methods Partial Differential Equations}, 22
(2006), 797-814. 

\bibitem{MP}
A.~ Mohamed and A.K.~Pani, 
An H1-Galerkin mixed finite element method combined with the modified method of characteristics for incompressible miscible displacement problems in porous media. New directions in applied mathematics (Hyderabad, 1995), 
{\em  Differential Equations Dynam. Systems}, 6(1998), 135-147. 

\bibitem{N}
L.~Nirenberg, An extended interpolation inequality,
{\em Ann. Scuola Norm. Sup. Pisa(3)}, 20 (1966), 733-737

\bibitem{Pea}
D.W.~Peaceman, {\em Fundamentals of Numerical Reservior Simulations}, Elsevier, Amsterdam,
1977.


\bibitem{RT}
P.A.~Raviart and J.M.~Thomas,
A mixed finite element method for 2nd order elliptic problems,
{\em Mathematical Aspects of Finite Element Methods, Lecture Notes in Math. 606,
Springer-Verlag, Berlin}, (1977), 292-315.

\bibitem{Rus}
 T. F. Russell,
Time stepping along characteristics
with incomplete iteration for a Galerkin approximation of miscible
displacement in porous media,
{\em SIAM J. Numer. Anal.}, 22 (1985), 970-1013.


\bibitem{SWML}
G. Scovazzi, M.F. Wheeler, A. Mikelic and S. Lee, 
Analytical and variational numerical methods for unstable miscible displacement flows in porous media, {\em J. Comput. Phys.}, 335(2017), 444--496.

\bibitem{SWS}
 Z. Si, J. Wang and W. Sun.
Unconditional stability and error estimates of
modified characteristics FEMs for the Navier-Stokes equations,
{\em Numer. Math.}, 134(2016), 139--161.

\bibitem{SY}
T.~Sun and Y.~Yuan, An approximation of incompressible miscible displacement in porous
media by mixed finite element method and characteristics-mixed finite element method,
{\em J. Comput. Appl. Math.}, 228 (2009), 391-411.

\bibitem{SW}
W.~Sun and C. Wu, 
New analysis and optimal error estimates of Galerkin-mixed FEMs for incompressible miscible flow in porous media, 
{\em Math. Comput.}, 2021
DOI: https://doi.org/10.1090/mcom/3561. 
	
\bibitem{Wang}
H.~Wang, An optimal-order error estimate for a family of ELLAM-MFEM approximations to
porous medium flow, {\em SIAM J. Numer. Anal.}, 46 (2008), 2133--2152.


\bibitem{WER}
H.~Wang, R.E.~Ewing and T.F.~Russell,
Eulerian-Lagrangian localized adjoint methods for convection-diffusion equations and their convergence analysis,
{\em IMA J. Numer. Anal.}, 15 (1995), 405--459.

\bibitem{WSS}
J. Wang, Z. Si and W. Sun, 
 A new error analysis of characteristics-mixed FEMs for miscible displacement in porous media, 
 {\em SIAM J. Numer. Anal.}, 52(2014), 3300--3020.
 
\bibitem{Whe}
M.F.~Wheeler, A priori $L^2$ error estimates for Galerkin approximations to parabolic partial
differential equations, {\em SIAM J. Numer. Anal.}, 10 (1973), 723--759.


\end{thebibliography}
\end{document}